\documentclass[notitlepage,11pt, a4paper]{article}
\usepackage[utf8]{inputenc}

\usepackage{amsfonts, amssymb, amsmath}
\usepackage{mathtools}

\usepackage{graphicx}
\usepackage{amsthm}
\usepackage{caption}
\usepackage{subcaption}

\usepackage[colorlinks=true, pdfstartview=FitV, linkcolor=blue, citecolor=blue, urlcolor=blue,pagebackref=false]{hyperref}
\usepackage{varioref}
\usepackage[capitalize]{cleveref}

\usepackage[numbers]{natbib}

\usepackage{amsthm}
\usepackage{thmtools}

\declaretheorem[numberwithin=section]{theorem}
\declaretheorem[sibling=theorem]{lemma}
\declaretheorem[sibling=theorem]{corollary}
\declaretheorem[sibling=theorem]{proposition}

\declaretheorem[sibling=theorem]{fact}

\usepackage{xcolor}
\definecolor{darkred}{rgb}{0.9,0.1,0.1}

\def \D {{\mathbb D}}
\def \E {{\mathbb E}}
\def \N {{\mathbb N}}
\def \P {{\mathbb P}}

\def \R {{\mathbb R}}

\def \cF {{\mathcal F}}

\def \cK {{\mathcal K}}
\def \cL {{\mathcal L}}

\def\one{\rlap{\mbox{\small\rm 1}}\kern.15em 1} 

\usepackage{stmaryrd} 

\title{Convergence of the height process of supercritical Galton--Watson forests with an application to the configuration model in the critical window}
\author{Serte Donderwinkel\thanks{\href{mailto:serte.donderwinkel@stats.ox.ac.uk}{serte.donderwinkel@stats.ox.ac.uk}, University of Oxford, Oxford, United Kingdom. ORCID iD: 0000-0001-8148-8631.}}
\date{\today}

\begin{document}
\maketitle
\begin{abstract}
\noindent We show joint convergence of the \L ukasiewicz path and height process for slightly supercritical Galton--Watson forests. This shows that the height processes for supercritical continuous state branching processes as constructed by Lambert (2002) are the limit under rescaling of their discrete counterparts. Unlike for {(sub-)critical} Galton--Watson forests, the height process does not encode the entire metric structure of a supercritical Galton--Watson forest. We demonstrate that this result is nonetheless useful, by applying it to the configuration model with an i.i.d.\ power-law degree sequence in the critical window, of which we obtain the metric space scaling limit in the product Gromov--Hausdorff--Prokhorov topology, which is of independent interest.
\end{abstract}
\section{Introduction}
In this work, we study the scaling limit of the  genealogical structure of a slightly supercritical Galton--Watson forest by showing convergence of its height process under rescaling. The height process is defined using a depth-first exploration of the forest. Consequently, it encodes only a part of a supercritical Galton--Watson forest: a walker that executes a depth-first exploration will never cross the first path to infinity that it encounters. We show that our convergence result can nevertheless be applied by using it to establish metric space convergence of the configuration model in the critical window. Moreover, our result shows that the height process for supercritical continuous state branching processes (CSBP) that was defined by Lambert in \cite{Lambert2002} is in fact the limit under rescaling of its discrete counterpart. 

We will give a short introduction to the encoding of a (random) forest by a \emph{\L ukasiewicz path} and a \emph{height process} in the discrete setting. Furthermore, we will introduce the \L ukasiewicz path and height process in the continuum. We then state our results and methods, followed by an overview of related work.

\subsection{Encoding forests by processes}\label{subsec.encodingrandomforests}
The encoding of random trees and forests in the discrete setting and the continuum by (excursions of) random processes has been around for a long time, see e.g.\ \cite{MeyerDwass1975,Rogers1984,NeveuPitman1989,Aldous1993,LeGall1998,Duquesne2005,Duquesne2007}. 

Let $T$ be an ordered rooted finite tree, say $|T|=n$. Let $v_0,\dots,v_{n-1}$ denote the vertices of the tree visited in depth-first order, so that $v_0$ is the root of the tree. \\
We will define the height function and \L ukasiewicz path of ${T}$. Both of these functions uniquely characterize ${T}$. The height process of ${T}$, referred to as $h$, is defined as $$h(i)=d_T(v_i,v_0),$$ i.e.\ for all $i$, $h(i)$ equals the distance from $v_i$ to the root.
Moreover, for all $i\in\{1,\dots,n\}$, let $y_i$ be the number of children of $v_{i-1}$, and set $y_0=1$. Then, the \L ukasiewicz path of ${T}$ is defined by $$x(i)=\sum\limits_{0\leq j\leq i} (y_j-1)$$ for $i=0,\dots,n$. Thus, $x(i)$ is the total number of younger siblings of $v_i$ and its ancestors, where \emph{younger} means that they are explored later in the depth-first search. Note that \begin{equation}\label{eq.heightprocessdef} h(i)=\#\{j< i:x(j)=\min\{x(k) : j\leq k \leq i\}\},\end{equation}
which is proved by Le Gall in \cite{LeGall1991}.
These encodings of trees by walks can easily be extended to ordered forests by concatenating the \L ukasiewicz paths and height processes.\\
We can use the correspondence between forests and their \L ukasiewicz paths to construct Galton--Watson forests from random walks. Suppose $(D_1,D_2,\dots)$ are i.i.d.\ random variables with law $\pi$ on $\N$. Then,
$$(S(k),k\geq 0):=\left(\sum_{i=1}^k(D_i-1),k\geq 0\right)$$
is the \L ukasiewicz path of a random forest in which all vertices have independent offspring with law $\pi$. We refer to such a forest as a \emph{$\pi$-Galton--Watson forest} ($\operatorname{GW}(\pi)$). We will write $(H(k),k\geq 0)$ to denote the height process corresponding to $(S(k),k\geq 0)$. \\

The continuous counterpart of a Galton--Watson process is a \emph{continuous-state branching process} (CSBP) (see e.g.\ \cite{GreyD1974}). The rôle of the \L ukasiewicz path is played by a Lévy process $(X_t,t\geq 0)$ without negative jumps (i.e.\ \emph{spectrally positive}). The law of $X$ is completely characterized by its Laplace exponent $\phi$, defined via 
$$\E\left[\exp({-\theta X_t})\right]=\exp(t\phi(\theta)).$$
We refer to the CSBP related to $X$ as $\operatorname{CSBP}(\phi)$. We wish to define the height process corresponding to $(X_t,t\geq 0)$ analogously to \eqref{eq.heightprocessdef}, so we should define a functional $H(X)=(H_t,t\geq 0)$ of $X$ such that $H_t$ in some sense measures the ``size'' of the set $\{s\leq t: X_{s-}=\inf\{X_r:r\in[s,t]\}$. In \cite{LeGall1998}, it was established that if $X$ almost surely does not drift to $\infty$ and satisfies
$$\int_1^\infty \frac{du}{\phi(u)}<\infty,$$ 
then there exists a continuous process $(H_t,t\geq 0)$ such that 
$$H_t{=}\lim_{\epsilon\to0}\frac{1}{\epsilon} \int_0^t\one_{\left\{X_s-\inf\{X_r:r\in[s,t]\}\leq \epsilon\right\}},$$
with the limit in probability for all $t\in[0,\infty)$. Further results were proved in \cite{Duquesne2005}. The excursions of $(H_t,t\geq 0)$ above $0$ encode metric spaces called \emph{Lévy trees}, which are defined in \cite{DuquesneLevytrees}. In \cite{Lambert2002}, the definition of $H(X)$ was extended to spectrally positive Lévy processes $X$ that drift to $\infty$ almost surely. 

\subsection{Results and methods}
For each $n$, let $D_1^n,D_1^2,\dots$ be an i.i.d.\ sequence of random variables with law $\pi_n$, and set $$(S^n(k),k\geq 0)=\left(\sum_{i=1}^k (D_i^n-1),k\geq 0\right).$$ 
Let $(H^n(k),k\geq 0)$ be the corresponding height process as defined in Subsection \ref{subsec.encodingrandomforests}. 
We impose the following conditions. 
\begin{description}
    \item[(C1)] 
    There exist a nondecreasing sequence of positive integers  $(\gamma_n,n\geq 1)$ that converges to $\infty$ and a L\'evy process $(X_t,t\geq 0)$ on $\R$ that does not have downward jumps and is of infinite variation, such that
    \begin{equation}\label{eq.assumptionconvergence}\left(n^{-1}S^n_{\lfloor n \gamma_n t \rfloor}, t\geq 0\right) \overset{d}{\to}(L_t,t\geq 0)\end{equation}
    as $n\to\infty$.
    \item[(C2)] For $Y^n_m$, the number of vertices at height $m$ in the first $n$ trees of the forest encoded by $(S^n(k),k\geq 0)$, for every $\delta>0$, \begin{equation}\label{eq.assumptiontreesdie}\liminf\limits_{n\to\infty}\P\left[Y^n_{\lfloor \delta \gamma_n\rfloor}=0\right]>0.\end{equation}
\end{description}
The following result is proved in \cite{Duquesne2005}.
\begin{theorem}\label{thm.convergencesubcritical}{\rm \cite[Theorem 1.4.3 and Theorem 2.2.1]{Duquesne2005}}
Suppose {\bf(C1)} and {\bf(C2)} hold, and in addition, we have
\begin{description}
\item[(C3*)] The Laplace exponent $\phi(\theta)$ of $(X_t, t\geq 0)$
satisfies 
\begin{equation*}\int_1^\infty \frac{du}{\phi(u)}<\infty.\end{equation*}
    \item[(C4*)] Almost surely, $X_t$ does not drift to $\infty$ as $t\to\infty$,
    \item[(C5*)] $\E[D^n_1]\leq 1$ for all $n$.
\end{description}
Then, there exists a continuous modification of the height process of $(X_t,t\geq 0)$, say $(H_t,t\geq 0)$, such that 
$$\left(\gamma_n^{-1}H^n_{\lfloor n \gamma_n t \rfloor}, t\geq 0\right)\overset{d}{\to}(H_t,t\geq 0)$$
as $n\to \infty$, jointly with \eqref{eq.assumptionconvergence}. 
\end{theorem}
Our main result is the analogue of Theorem \ref{thm.convergencesubcritical} for supercritical Galton--Watson processes.
\begin{theorem} \label{thm.mainconvergenceresult}
 If {\bf(C1)} and {\bf(C2)} hold, and additionally,
\begin{description}
\item[(C3)] 
For $\phi(\theta)$ the Laplace exponent of $(X_t, t\geq 0)$ and $\xi>0$ the unique value such that $\phi(\xi)=0$, we have 
\begin{equation*}\int_1^\infty \frac{du}{\phi(u+\xi)}<\infty.\end{equation*}
    \item[(C4)] $X_t\to \infty$ almost surely as $t\to\infty$,
    \item[(C5)] $\E[D^n_1]>1$ for all $n$.
\end{description}
Then, there exists a continuous modification of the height process of $(X_t,t\geq 0)$, say $(H_t,t\geq 0)$, such that 
$$\left(\gamma_n^{-1}H^n_{\lfloor n \gamma_n t \rfloor}, t\geq 0\right)\overset{d}{\to}(H_t,t\geq 0)$$
as $n\to \infty$ jointly with \eqref{eq.assumptionconvergence}. 
\end{theorem}
Roughly speaking, condition {(C2)} ensures that the maximal height in a forest of a number $n$ of $\operatorname{GW}(\pi_n)$ trees, rescaled by $\gamma_n^{-1}$, may not exceed $\delta$. An analogous fact holds for $\operatorname{CSBP}(\phi)$.  Informally, condition {(C3)} ensures that a Lévy tree encoded by an excursion of $X_t$ above its infimum, conditional on the excursion being finite, has finite height almost surely. This is a necessary condition for a continuous modification of the height process to exist. Finally, conditions {(C4)} and {(C5)} ensure that $\operatorname{GW}(\pi_n)$ and $\operatorname{CSBP}(\phi)$ respectively are supercritical.

Our method is based on the following pathwise construction of $(S^n(k),H^n(k), k\geq 0)$. 
\begin{enumerate}
    \item Up to the first time that it hits it overall infimum, $(S^n(k),k\geq 0)$ encodes a random number of finite trees, which are independent, and each have the law of a tree in $\operatorname{GW}(\pi_n)$ conditioned to be finite. This yields a pathwise construction of $(S^n(k),H^n(k), k\geq 0)$ up to the first time that $(S^n(k),k\geq 0)$ hits its overall infimum. 
    \item After the first time that it hits its overall infimum, $(S^n(k),k\geq 0)$  encodes an infinite spine with trees attached to the left of the infinite spine, and vertices attached to the right of the infinite spine. To be precise, for every vertex $v$ on the infinite spine, a random number of independent copies of a tree in $\operatorname{GW}(\pi_n)$ conditioned to be finite are attached to $v$ left of the infinite spine, and a random number of vertices (that will never be visited) attached to the right of the infinite spine. This yields a pathwise construction of $(S^n(k),H^n(k), k\geq 0)$ after the first time that $(S^n(k),k\geq 0)$ hits its overall infimum. This pathwise construction is similar to the pathwise construction for Galton--Watson processes with immigration defined in \cite{Duquesne2007}. 
\end{enumerate}
(See e.g.\ \cite{Janson2012}, \cite[Section 5.7]{Lyons2017}, and \cite[Chapter 12]{Athreya1972}  for the laws of the trees encoded by $(S^n(k), k\geq 0)$ before and after it hits its overall infimum, although height processes are not considered in these works.)

We define a similar pathwise construction of $(L_t,H_t,t\geq 0)$, which is standard for the pre-infimum process (see e.g.\ \cite{Bertoin1996}, Chapter VII), and based on  \cite{Lambert2002}, Section 5 for the post-infimum process. We then show convergence in distribution under rescaling of the pathwise construction of $(S^n(k),H^n(k), k\geq 0)$ to the pathwise construction of $(L_t,H_t,t\geq 0)$. 

Finally, in Section \ref{sec.application} we use Theorem \ref{thm.mainconvergenceresult} to show metric space convergence of the largest components of a uniform graph with an i.i.d.\ heavy-tailed degree sequence in the critical window, which extends the main result of \cite{Conchon2018}. For a precise statement of this result, and an overview of earlier work on related graph models, we refer the reader to Section \ref{sec.application}. 

\subsection{Related work}
While there is an extensive literature concerning the convergence of height processes for critical branching processes, we are only aware of two works that consider the supercritical case: \cite{broutinLimitsMultiplicativeInhomogeneous2020} and \cite{Duquesne2009}. In  \cite[Theorem B2]{broutinLimitsMultiplicativeInhomogeneous2020}, Broutin, Duquesne and Wang discuss convergence of the height process of a model similar to the model we consider here. However, in that work, the authors only consider the forest before the first infinite line of descent. In \cite{Duquesne2009}, Duquesne studies a class of supercritical branching processes with an single infinite line of descent (\emph{sin-trees}). The author encodes the resulting trees with two processes, one that encodes the genealogical structure on the left, and one that encodes the genealogical structure on the right side of the infinite line of descent. For this model, convergence under rescaling of the discrete contour processes to the continuum height process is proved. The pathwise construction that we use resembles the pathwise construction used in \cite{Duquesne2009}. 

Convergence of supercritical Galton--Watson trees under rescaling was studied through a different lens in \cite{Duquesne2012}. In Theorem 4.15 of that work, Duquesne and Winkel show the convergence of a class of supercritical Galton--Watson forests to a Lévy forest in the sense of Gromov--Hausdorff convergence on locally compact rooted real trees. Although their theorem applies to the entire forest, and not just to the trees and pendant subtrees to the left of the first infinite line of descent, convergence in the Gromov--Hausdorff topology on locally compact rooted real trees does not imply convergence of ``depth-first ordering'' of the vertices in the tree, as convergence of the height process does.

An alternative approach to view the height process of a supercritical branching process was introduced in \cite{Delmas2008} for branching processes with a quadratic branching mechanism, and extended to general CSBPs in \cite{Abraham2012pruning}. The approach is to build the super-critical tree up to a given level $a$, such that the tree can be encoded by a height process. Then, the law of the resulting tree is absolutely continuous with respect to the law of a (sub-)critical Lévy tree whose branches above level $a$ are removed, which is referred to as \emph{pruning}. Furthermore, this family of processes indexed by $a$ satisfies a consistency property, and hence there exists a projective limit. It is established in \cite{Abraham2015} that the limit object has the same law as the supercritical Lévy tree defined in \cite{Duquesne2007}. In \cite{HeHui2013}, He and Luan define an analogous pruning operation for supercritical Galton--Watson trees and prove that the contour functions of these truncated Galton--Watson trees converge weakly to the processes constructed by Abraham and Delmas.

\section{The construction of the height process}\label{sec.construction}

In this section, we will combine results from the literature in order to give a construction of the height process of a supercritical Galton--Watson process in the discrete case, and of a supercritical CSBP in the continuum.

\subsection{The height process in the discrete case}\label{subsec.constdiscrete}
In this section, we will describe a pathwise construction of the height process corresponding to $S^n$. We do this by considering the process before and after it hits its overall infimum separately. This corresponds to considering the laws of finite and infinite trees in a supercritical branching process separately; this idea can be traced back to Harris \cite{Harris1948}. In the descriptions of the two parts, we will make use of a process that is locally absolutely continuous to $S^n$, which we will denote by $\hat{S}^n$. We will refer to to $\hat{S}^n$ as ``$S^n$ conditioned to die out". This formalizes, in the random walk framework, the well-known fact that a supercritical Galton--Watson process conditioned to die out is a subcritical Galton--Watson process. \\
For $t\geq 0$, define 
$$\cL_{S^n(1)}(\theta)=\E\left[\exp(-\theta S^n(1))\right],$$
and set $\phi_n(\theta)=\log \cL_{S^n(1)}(\theta)$. 
Since $\E[S^n(1)]>0$, we have $\cL_{S^n(1)}'(0)<0$, so by the convexity of $\cL_{S^n(1)}(\theta)$ and the fact that $\cL_{S^n(1)}(0)=1$, there is a unique $\xi_n>0$ such that $\cL_{S^n(1)}(\xi_n)=1$ and  $\phi_n(\xi_n)=0$. Let $\P_n$ be the law of $S^n$, and let 
$$\cF_k^n:=\sigma(S^n(m),m\leq k)$$
be the natural filtration of $S^n$. Let $\P_n^\#$ be locally absolutely continuous with respect to $\P_n$, with
$$\mathrm{d}\P_n^\#|_{\cF_k^n}=\exp\left(-\xi_n S^n(k)\right) \mathrm{d}\P_n|_{\cF_k^n},$$
and let $\hat{S}^n$ be a process which under $\P_n$ has the law of $S^n$ under $\P_n^\#$. Let $$\tau^n(m)=\inf\left\{k:S^n(k)=-m\right\},$$ so that $\{\tau^n(m)<\infty\}$ may be interpreted as the event that at least $m$ trees in the Galton--Watson forest are finite. The following properties of $\xi_n$ and $\hat{S}^n$ are standard (see e.g.\ \cite{Athreya1972}, Chapter 12). 
\begin{theorem}\label{thm.factsonprocessdyingout} We have
    \begin{enumerate}
        \item for any $m\geq 0$, \begin{align*} \exp(-\xi_n m)&=\P_n\left[\tau^n(m)<\infty \right];
    ,\end{align*}
   \item  if $\Lambda\in \cF_{\tau^n(m)}^n$ for some $m>0$, then
    \begin{align*} \P^\#_n[\Lambda]=\P_n\left[\Lambda|\tau^n(m)<\infty\right]; \text{ and}\end{align*}
    \item $\hat{S}^n$ is a downward skip-free random walk on the integers.
    \end{enumerate}
\end{theorem}
The following lemma is a first reason why $\hat{S}^n$ plays an important rôle in the pathwise construction of $S^n$. 

\begin{lemma}\label{lemma.preinfdiscrete}
Let $G^n$ be a geometric random variable with success probability $\exp(-\xi_n)$, i.e. $\P_n(G^n=k)=\exp(-k\xi_n)(1-\exp(-\xi_n))$, independent of $\hat{S}^n$. Then, the pre-infimum process of $S^n$ has the law of $\hat{S}^n$, stopped at the first time it reaches level $G^n$.
\end{lemma}
\begin{proof}
Note that the negative of the overall infimum of $S^n$, which we denote by $-I^{n}_{\infty}$, equals the number of finite trees in the forest defined by $S^n$ viewed as a \L ukasiewicz path. Hence, it is distributed as a geometric random variable with parameter $\exp(-\xi_{n})$.  Let $\rho^n$ denote the time when $S^n$ first reaches $I^{n}_{\infty}$. We have that $$(S^n(j):0\leq j \leq \rho^n)\text{ under }\P_{n}(\;\cdot\;|I^n_{\infty}=-m) $$ has the same distribution as $$(S^n(j): 0\leq j \leq \tau^n(m))\text{ under }\P_{n}\left(\;\cdot\;|\tau^n(m)<\infty\right),$$ which, by Theorem \ref{thm.factsonprocessdyingout},  is equal in distribution to $$(\hat{S}^{n}(j):0\leq j \leq \tau^n(m))\text{ under }\P_{n}^\#.$$ Combining this with the distribution of $I^{n}_{\infty}$, we find that for $G_n$ a random variable with the geometric distribution with success probability $\exp(-\xi_n)$ independent of $S^n$, $(S^n(j):0\leq j \leq \rho^{n})$ under $\P_{n}$ has the same distribution as $(\hat{S}^{n}(j):0\leq j \leq \tau^n(G^n))$ under $\P_{n}^\#$. 
\end{proof}

By Theorem \ref{thm.factsonprocessdyingout},  $\hat{S}^n$ encodes a subcritical Galton--Watson forest. We let $\hat{H}^n$ be its height process. Then, Lemma \ref{lemma.preinfdiscrete} has the following corollary.
\begin{corollary}\label{cor.preinf}
We have that 
$$(S^n(m),H^n(m),m\leq \rho^n)\overset{d}{=}(\hat{S}^n(m),\hat{H}^n(m),m\leq \tau^n(G^n)).$$
\end{corollary}

We will now characterise the post-infimum process and its height process. Let $\P_n^\uparrow$ be the law of $(S^n(k),k\geq 0)$ conditioned to be non-negative for all $k$. Note that we are conditioning on an event of non-zero probability because $S^n$ is supercritical, so this is well-defined. 

The following lemma suggests an adaptation of $S^n$ that has the same height process, and that turns out to be easier to work with.

\begin{lemma}\label{lemma.heightprocessesagree} For a discrete \L ukasiewicz path $(X_k,k\geq 0)$ set $\underline{\underline{X}}_n=\min\{X_m, m\geq n\}$. We will refer to $\underline{\underline{X}}_n$ as the \emph{future infimum process}. Then, the height processes of $X$ and $X-\underline{\underline{X}}$ are the same. \end{lemma}
\begin{proof}
Let  $H^*$ denote the height process of $X-\underline{\underline{X}}$. Then, by definition, 
$$H^*_n=\#\left\{m< n : X_m-\underline{\underline{X}}_m=\min\{X_j-\underline{\underline{X}}_j , m\leq j \leq n \}\right\}$$
Set $\underline{\underline{g}}(n)=\max\{ m < n : X_m=\underline{\underline{X}}_m\}$. Then, for $m\leq \underline{\underline{g}}(n)$,  we have that $\underline{\underline{X}}_m=\min\{ X_j: m\leq j \leq n\}$ and $\min\{X_j-\underline{\underline{X}}_j , m\leq j \leq n \}=0$. Also, for $m>\underline{\underline{g}}(n)$, we have $\underline{\underline{X}}_m=\underline{\underline{X}}_n.$
Hence, indeed,
\begin{align*}
    H^*_n=&\#\left\{m\leq \underline{\underline{g}}(n) : X_m-\underline{\underline{X}}_m=\min\{X_j-\underline{\underline{X}}_j , m\leq j \leq n \}\right\}\\
    &+\#\left\{\underline{\underline{g}}(n)<m<n : X_m-\underline{\underline{X}}_m=\min\{X_j-\underline{\underline{X}}_j , m\leq j \leq n \}\right\}\\
    =&\#\left\{m\leq \underline{\underline{g}}(n) : X_m-\min\{ X_j: m\leq j \leq n\}=0\right\}\\&+\#\left\{\underline{\underline{g}}(n)<m<n : X_m-\underline{\underline{X}}_n=\min\{X_j-\underline{\underline{X}}_n , m\leq j \leq n \}\right\}\\
    =&\#\left\{m\leq n : X_m=\min\{X_j , m\leq j \leq n \}\right\} = H_n,
\end{align*}
where $H_n$ is the height process of $X_n$. \end{proof}

\begin{figure}
\centering
\includegraphics{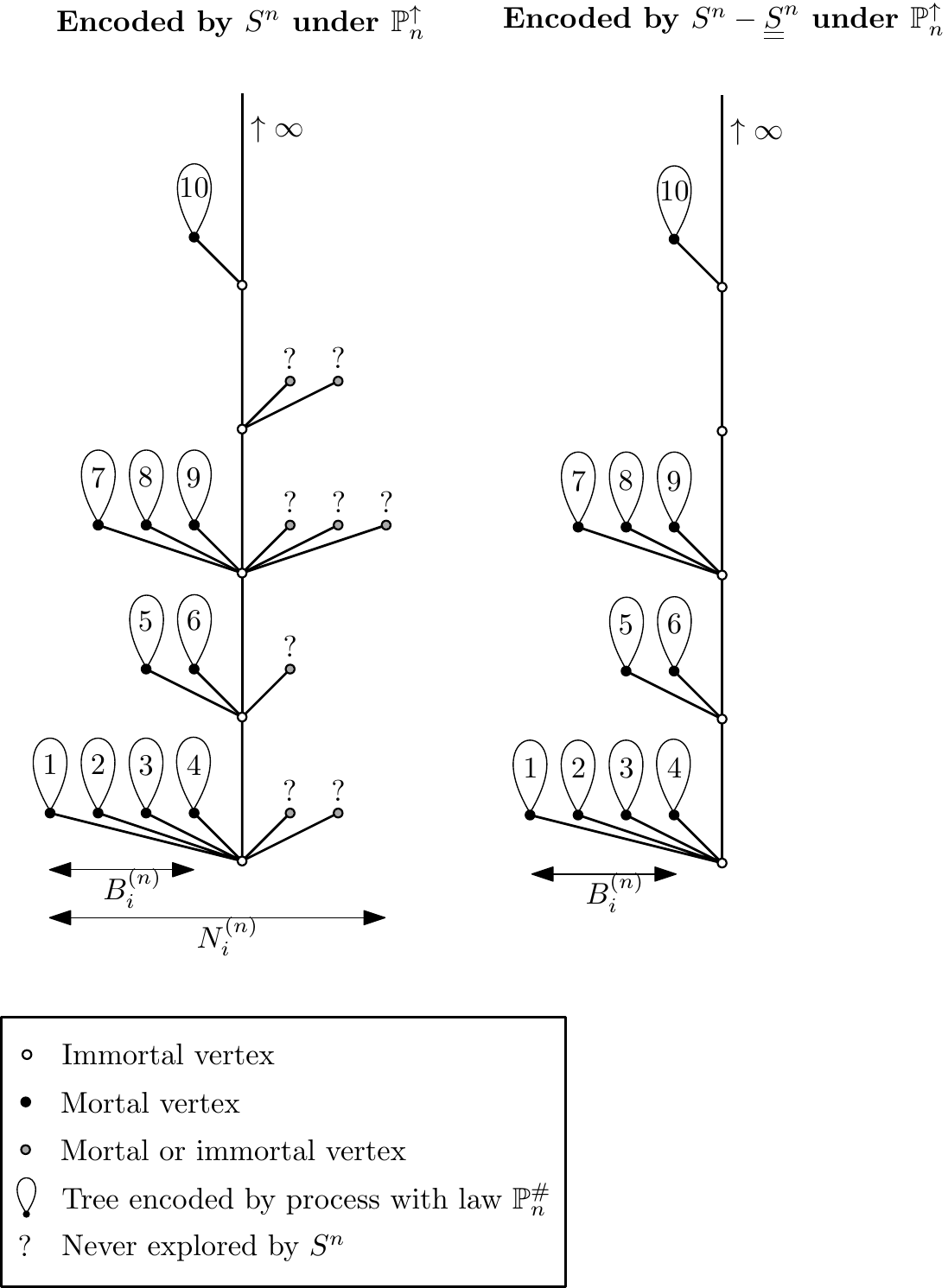}
\caption{The information captured in $S^n$ and $S^n-\underline{\underline{S}}^{n}$ under $\mathbb{P}_n^\uparrow$.}
\label{figure.trees}
\end{figure}

Note that since $S^n$ drifts to $+\infty$ almost surely, $S^n-\underline{\underline{S}}^{n}$ under $\P^\uparrow_{n}$ is recurrent in zero. Indeed, if $S^n$ drifts to $+\infty$, we have  $\underline{\underline{S}}^{n}(k)>-\infty$ for all $k>0$. Then, by definition of $\underline{\underline{S}}^{n}$, there is a finite $l>k$ such that $S^n(l)=\underline{\underline{S}}^{n}(l)$. Moreover, $S^n-\underline{\underline{S}}^{n}$ is non-negative. Combining these facts implies that $S^n-\underline{\underline{S}}^{n}$ is the \L ukasiewicz path of an infinite spine with finite trees attached to it only on the left-hand side. $S^n$ describes the same infinite spine with trees on its left-hand side, but also contains information on the total number of children (including children to the right of the spine) of vertices on the spine. See Figure \ref{figure.trees}. We will use this description to study the distributions of these paths.

We will use terminology introduced by Janson in \cite{Janson2012} and will call a vertex \emph{immortal} if it is the root of an infinite tree. Otherwise we will call it mortal. Note that since $S^n$ drifts to $+\infty$ almost surely, every vertex has a positive probability of being immortal. Hence, there are almost surely countably infinitely many infinite spines. Note that in every generation, we only visit the vertices to the left of the leftmost immortal vertex. This means that, under $\P^\uparrow_{n}$, an excursion of $S^n-\underline{\underline{S}}^{n}$ above zero consists of an increment (of size, say, $B^{(n)}$) corresponding to the mortal older brothers of the first immortal vertex, and then an excursion with the law of $S^n$ starting at $B^{(n)}$ conditioned to hit $0$ in finite time. This corresponds to first sampling the number of trees to the left of the infinite spine and then the shapes of those trees. \\
The sample paths of $S^n$ can be constructed from the sample paths of $S^n-\underline{\underline{S}}^{n}$, by adding the randomness that encodes the number of younger brothers of vertices on the leftmost path to infinity. This corresponds to replacing the jumps of size $B^{(n)}$ by jumps of size $N^{(n)}-1$, with $N^{(n)}$ being distributed as the total generation size given $B^{(n)}$. This is illustrated in Figure \ref{figure.trees}.
(A similar decomposition of Galton--Watson trees conditioned on non-extinction is discussed by Lyons and Peres in Section 5.7 of \cite{Lyons2017}, where they do not introduce the encoding by a \L ukasiewicz path and height process.)
\\
Note that the height at time $k$ in the tree encoded by $S^n$ under $\P^\uparrow_n$ is given by the height of the position along the infinite spine where the finite subtree containing the vertex visited at time $k$ is attached plus the height of the vertex in the finite subtree.\\
 We will now investigate the joint distribution of $N^{(n)}$ and $B^{(n)}$.
 \begin{lemma}\label{lemma.brotherdistr}
 Let $N^{(n)}$ be the number of children in a set of offspring conditioned to contain at least one immortal vertex, let $B^{(n)}$ be the number of older brothers of the oldest immortal vertex in such a set of offspring, and let $\exp(-\xi_n)$ be the probability that, under $\P_{n}$, a tree dies out. Then, we have for $l>k$, 
\begin{equation}
 \P(B^{(n)}=k,N^{(n)}=l)=\exp(-k\xi_n)\P(S_1^n=l-1)
\end{equation}
 \end{lemma}
 \begin{proof}

 Let $\psi_n$ denote the probability generating function of the offspring distribution under $\P_{n}$ (i.e. the law of $Z^n_1$ under $\P_{n}$). Let $M_n$ be the random variable representing the number of mortal children of an immortal parent, and let $J_n$ be the random variable representing the number of immortal children of an immortal parent. Then,  
 in  \cite{Janson2012}, Janson gives the generating function of the joint law of $M_n$ and $J_n$ as
 \begin{equation}\label{eq.genfunjanson} \E\left[x^{M_n} y^{J_n}\right]=\frac{\psi_n\left(\exp(-\xi_n) x+(1-\exp(-\xi_n))y\right)-\psi_n\left(\exp(-\xi_n) x\right)}{1-\exp(-\xi_n)}.\end{equation} Also, given the number of mortal and immortal children, they appear in a uniformly random order.
It is then straightforward that the generating function of the total number of children of an immortal parent is given by 
 $$\E\left[x^{N^{(n)}}\right]=\frac{\psi_n(x)-\psi_n(\exp(-\xi_n) x)}{1-\exp(-\xi_n)},$$
 so we obtain that for $k=1,2,\dots$, \begin{equation*}\P(N^{(n)}=k)=\frac{1-\exp(-k\xi_n)}{1-\exp(-\xi_n)}\P(S_1^n=k-1).\end{equation*} Using \eqref{eq.genfunjanson} to analyse the generating function of the joint law of $N^{(n)}$ and $J_n$, we see that, conditional on the value of $N^{(n)}$, $J_n$ is distributed as a binomial random variable with parameters $N^{(n)}$ and $1-\exp(-\xi_n)$ conditioned to be at least $1$. Since the mortal and immortal children appear in a uniform order, conditional on the generation size $N^{(n)}$, the number of mortal older brothers of the first immortal vertex, $B^{(n)}$, is distributed as a geometric random variable with success parameter $\exp(-\xi_n)$ conditioned to be at most $N^{(n)}-1$. 
 We obtain that 
 $$\P(B^{(n)}=k|N^{(n)}=l)=\frac{\exp(-k\xi_n)(1-\exp(-\xi_n))}{1-\exp(-l\xi_n)},$$
 which proves the statement.
\end{proof}
These findings on the post-infimum process can be summarised as follows.
\begin{fact}\label{fact.lawdiscretepostinf}
The sample paths of $S^n$ and $S^n-\underline{\underline{S}}^n$ under $\P^\uparrow_n$ can be constructed by concatenating excursions above the future infimum as follows. Sample a countably infinite number of independent copies of $B^{(n)}$ and $N^{(n)}$ according to $\P(B^{(n)}=k,N^{(n)}=l)=\exp(-k\xi_n)\P(S_1^n=l-1)$ for $k<l$. Start an excursion with an increment of size $N^{(n)}-1$ above the previous future infimum, and continue from there as a process with law $\P^\#_n$. Kill the excursion at $N^{(n)}-B^{(n)}-1$ above the previous future infimum, which will be the new future infimum.\\
By replacing the jumps of size $N^{(n)}-1$ by jumps of size $B^{(n)}$ we obtain  $S^n-\underline{\underline{S}}^{n}$. 
\end{fact}
\subsubsection{A pathwise construction}\label{subsubsec.pathwiseconstr}
The characterisation of the pre- and post-infimum process justifies the following pathwise construction of $(S^n,S^n-\underline{\underline{S}}^n, H^n)$ under $\P_n$. This is similar to the pathwise construction of the encoding processes of a Galton-
Watson process with immigration given in Section 2.2 of \cite{Duquesne2009}.

\begin{figure}
\centering
\includegraphics{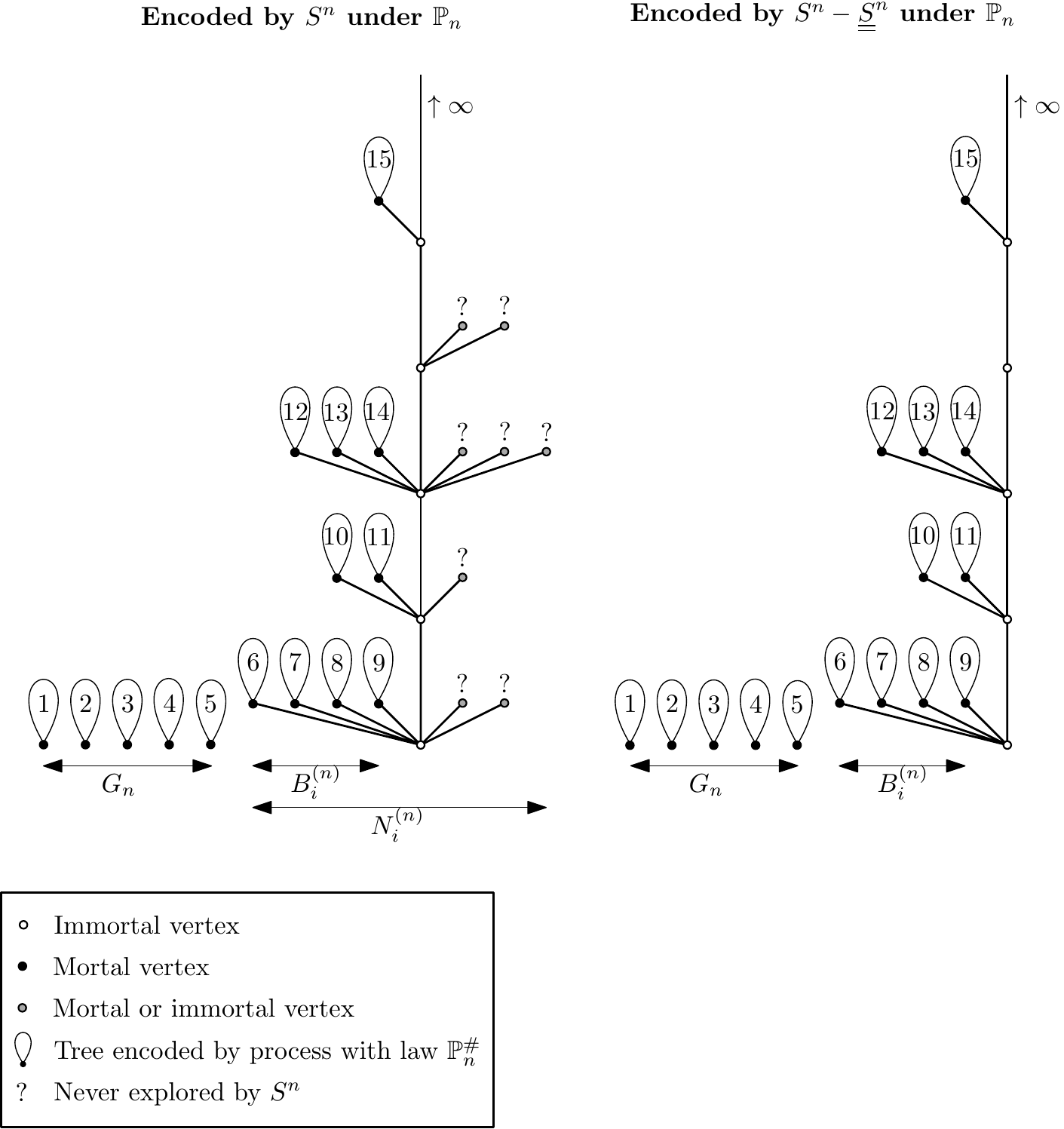}
\caption{The information captured in $S^n$ and $S^n-\underline{\underline{S}}^{n}$ under $\mathbb{P}_n$. The finite trees are encoded by the pre-infimum process, and the infinite spine and its pendant subtrees are encoded by the post-infimum process.}
\label{figure.preandpostcombined}
\end{figure}

See Figure \ref{figure.preandpostcombined} for a graphical representation of the construction.

\begin{itemize}
    \item Let $\hat{S}^n(k)$ be a random walk with law $\P^\#_n$ and let $\hat{H}^n(k)$ be the corresponding height process. Set $\hat{I}^n(k)=\min_{m\leq k} \hat{S}^n(k)$. The trees encoded by $\hat{S}^n(k)$ will be used as the finite trees that are explored before the first infinite tree is encountered, and as the finite subtrees to the left of the left-most infinite spine which are rooted at the vertices on the infinite spine.
   \item Let $G^n$ be distributed as a geometric random variable with mean $\exp(-\xi_n)$, which  gives the number of finite trees explored by the process. Let $(B^{(n)}_1,N^{(n)}_1),(B_2^{(n)},N^{(n)}_2),\dots$ be i.i.d.\ pairs of random variables, independent of $\hat{S}^n$, distributed as $(B^{(n)},N^{(n)})$. Set $Q^n(0)=G_n$, $Q^n(k)=G_n+\sum\limits_{i=1}^k N_i^{(n)}$ for $k\geq 1$. $N^{(n)}_m$ will be the number of children attached to the $m^\text{th}$ vertex on the infinite spine. Also, set $D^n(0)=G_n$ and $D^n(k)=G_n+\sum\limits_{i=1}^k B_i^{(n)}$ for $k\geq 1$. $B^{(n)}_m$ will be the number of finite subtrees to the left of the infinite spine rooted at the $m^\text{th}$ vertex on the infinite spine. 
  \item Define $F^{n}(k)=\inf\{l:-\hat{I}^n(k-1-l)<D^n(l)\}.$ $F^{n}(k)$ will be the number of vertices located on the leftmost infinite spine among the first $k$ vertices that we visit in our depth-first exploration.
   
    \item Then, define
        $$S^n(k)=-G^{(n)}+\hat{S}^n\left(k-{F}^{(n)}(k)\right)+{Q}^n\left({F}^{(n)}(k)\right)-F^{n}(k)$$
        so that
        $$(S^n-\underline{\underline{S}}^{n})(k)= \hat{S}^n(k-{F}^{(n)}(k))+{D}^n\left({F}^{(n)}(k)\right)$$
        and
        \begin{align}\begin{split}\label{eq.pathwiseconstrSandHentireprocess}H^n(k)&= \left(\hat{H}^n(k-{F}^{(n)}(k))+{F}^{(n)}(k)\right)\one_{\{{F}^{(n)}(k-1)={F}^{(n)}(k)\}}\\&+  \left(\hat{H}^n\left(k-({F}^{(n)}(k)-1)\right)+{F}^{(n)}(k)-1\right)\one_{\{{F}^{(n)}(k-1)<{F}^{(n)}(k)\}}.
        \end{split}
        \end{align}
        
\end{itemize}

By Corollary \ref{cor.preinf}, Fact \ref{fact.lawdiscretepostinf}, and Lemma \ref{lemma.heightprocessesagree}, the constructed process has the desired distribution.

In the next section, we will introduce the continuous counterpart of this construction. In Section \ref{sec.convergence}, we will then use the form \eqref{eq.pathwiseconstrSandHentireprocess} to show the joint convergence under rescaling of the discrete processes to the continuous processes.

\subsection{The height process of a supercritical L\'evy process}\label{sec.supercriticalheightprocess}

Just as in the discrete case, we will obtain a pathwise construction of a supercritical spectrally positive L\'evy process and its height process by considering its pre- and post-infimum processes separately. As before, let $L$ be such a L\'evy process with Laplace exponent $\phi$. Define $I_\infty=\inf\{L_t,t\geq 0\}$, and set $\rho=\sup\{t>0:L_{t}=I_\infty\}$. The process on $[0,\rho]$ will be referred to as the \emph{pre-infimum process}, and the process on $[\rho,\infty)$ will be referred to as the \emph{post-infimum process}. Informally, as in the discrete case, the pre-infimum process encodes the $\R$-trees without a path of infinite length, and the post-infimum process encodes the metric structure to the left of the leftmost path of infinite length in the first $\R$-tree that contains such a path. \\
(In \cite{Bertoin1993}, Bertoin uses a different strategy, and splits the process at the infimum attained on a compact time interval. We will not use this approach, and so we will not discuss his results here.) \\
On an infinite time horizon, the following results on spectrally positive L\'evy processes are available:
\begin{enumerate}
  \item By Th\'eor\`eme 2 in Bertoin \cite{Bertoin1991}, the pre-infimum process of a spectrally positive L\'evy process that drifts to $+\infty$ has the same law as the process `conditioned to drift to $-\infty$' stopped at an independent exponential level. Informally, this result says that the pre-infimum process encodes $\R$-trees conditioned to die out. 
    \item By Lemma 4.1 in  Millar \cite{Millar1977} (which is rephrased to the result we need as `Th\'eor\`eme (Millar)' in \cite{Bertoin1991}), the post-infimum process of a spectrally positive L\'evy process that drifts to $+\infty$ has the same law as the process conditioned to stay positive, and is independent of the pre-infimum process. Informally, this result says that the post-infimum process encodes (part of) a single $\R$-tree conditioned to be infinite.
    \item In \cite{Lambert2002}, Lambert describes the height process corresponding to a class of spectrally positive L\'evy processes conditioned to stay positive.
\end{enumerate}

We start with an overview of the results by Bertoin \cite{Bertoin1991} that we use. Firstly, since $L$ is supercritical, there is a unique $\xi>0$ such that $\phi(\xi)=0$. Let $\P$ be the law of $L$, and set 
$$\cF_t:=\sigma(L_s, s\leq t).$$ Then, $(\exp(\xi L_t),t\geq 0)$ is a $\P$-martingale. Let $\P^\#$ be locally absolutely continuous with respect to $\P$, with 
\begin{equation}\label{def.Phash}\mathrm{d}\P^\#|_{\cF_t}=\exp(-\xi L_t)\mathrm{d}\P|_{\cF_t}.\end{equation}
Let $\hat{L}$ be a process which under $\P$ has the law of $L$ under $\P^\#$. The following analogue of Theorem \ref{thm.factsonprocessdyingout} is a straightforward consequence of the fact that $\exp(\xi L_t)$ is a martingale. See \cite[Chapter VII]{Bertoin1996}.
\begin{theorem}\label{thm.propertiesconditionedprocess}
We have
\begin{enumerate}
    \item for $\tau(x)=\inf\{t>0 : L_t=-x\}$, $$\P[\tau(x)<\infty]=\exp(-\xi x),$$ 
    \item if $\Lambda\in \cF_{\tau(x)}$ for some $x>0$, then  
    $$\P^\#[\Lambda]=\P[\Lambda|\tau(x)<\infty]; \text{ and}$$
    \item $\P^\#$ is the law of a spectrally positive subcritical L\'evy process with Laplace exponent $\phi^\#(\cdot)=\phi(\cdot+\xi)$.
\end{enumerate}
\end{theorem}
The following theorem is then proved in \cite{Bertoin1991} as Th\'eor\`eme 2. 

\begin{theorem}[Th\'eor\`eme 2, \cite{Bertoin1991}]\label{thm.preinfcontlaw}
Let $E$ be an exponential random variable with rate $\xi$. The pre-infimum process of $L$ has the law of $\hat{L}$, independent of $E$, stopped at the first time it reaches level $E$.
\end{theorem}

These observations, together with Proposition 1.4.3 in Duquesne and Le Gall in \cite{Duquesne2005}, imply the following proposition.

\begin{proposition}\label{prop.preinfcont}
There exists a continuous modification of the height process corresponding to $\hat{L}$, which we will denote by $\hat{H}$. Moreover, for $$\rho=\sup\{t:s<t<u\implies L_t\leq L_s\text{ and } L_t\leq L_u\},$$
we have $\rho<\infty$ almost surely, and there exists a continuous modification of the height process of $L$ up to $\rho$, which we will refer to as $H$, and
$$(L_t,H_t,t\leq \rho)\overset{d}{=}(\hat{L}_t,\hat{H}_t, t\leq \tau(E)).$$
\end{proposition}
\begin{proof}
Since $L$ is supercritical, $\rho$ is finite almost surely.
Theorem \ref{thm.propertiesconditionedprocess}.3 and condition {(C3)} ensure that the conditions of Theorem 1.4.3 in \cite{Duquesne2005} are satisfied, which implies that the height process corresponding to $\hat{L}$ exists and has a continuous modification. Theorem  \ref{thm.preinfcontlaw} then yields the proposition.
\end{proof}

We will now focus on the post-infimum process. We will use two important results from the literature. Firstly, by Lemma 4.1 in Millar \cite{Millar1977}, the post-infimum process $(L_{\rho+t}-I_{\infty}, t\geq 0)$ has the same law as $L$ conditioned to stay positive and is independent of the pre-infimum process. Call the law of $L$ conditioned to stay positive $\P^\uparrow$. For the definition of this process, see \cite{Bertoin1991}, \cite{Chaumont1994}, and \cite{Chaumont1996}.  The height process of $L$ under $\P^\uparrow$ is characterised by Lambert in \cite{Lambert2002}, and is obtained via the continuous counterpart of the construction in Section \ref{subsec.constdiscrete}. Indeed, in \cite{Lambert2002}, Lambert defines $$\underline{\underline{L}}_t:=\inf\{L_s,s\geq t\},$$ and, in Lemma 8ii \cite{Lambert2002}, he shows that the height processes of $L-\underline{\underline{L}}$ and $L$ are equal. Then, in Theorem 3 \cite{Lambert2002}, he gives a pathwise construction of $L-\underline{\underline{L}}$ and its local time at $0$ under $\P^\uparrow$, by viewing $L-\underline{\underline{L}}$ as a continuous time branching process with immigration, an object introduced in \cite{Kawazu1971}. Lemma 8ii \cite{Lambert2002} also illustrates how to construct the height process corresponding to $L-\underline{\underline{L}}$.  Translating these results to our setting yields the following proposition, which is the continuous counterpart of the pathwise construction of $S^n-\underline{\underline{S}}^n$ under $\P_n^\uparrow$  discussed in Fact \ref{fact.lawdiscretepostinf}. 

\begin{theorem}[Theorem 3, Lemma 8ii  \cite{Lambert2002}] \label{thm.postinfcont}
Recall the definition of $\P^\#$ from \eqref{def.Phash}, and let $\hat{L}$ be a process which under $\P$ has the law of $L$ under $\P^\#$, and let $\phi^\#$ be its Laplace exponent. Let $\widetilde{D}$ be a spectrally positive L\'evy process with Laplace exponent $\frac{\phi(\cdot)}{\cdot+\xi}$ independent of $\hat{L}$. Then, $\widetilde{D}$ is a subordinator. For $t\geq 0$, define the right-inverse of $\widetilde{D}$
$$\widetilde{D}^{-1}_t:=\inf\{s:\widetilde{D}_s>t\}$$
and set 
$$\hat{I}_t=\inf_{s\leq t}\hat{L}_s.$$
Then, for $L^*$ defined by 
$$L^*_t=\hat{L}_t+\widetilde{D}\left(\widetilde{D}^{-1}(\hat{I}_t)\right),$$
$L-\underline{\underline{L}}$ under $\P^\uparrow$ has the same law as $L^*$.
Moreover, the local time of $L^*$ at zero may be defined by 
$$\ell^*_t=\widetilde{D}^{-1}(\hat{I}_t).$$
Finally, suppose that $\hat{H}$ is a continuous modification of the height process corresponding to $\hat{L}$. Then, 
$$H^*_t:=\ell^*_t+\hat{H}_t$$
is a continuous modification of the height process corresponding to $L^*_t$.

\end{theorem}
Combining the above proposition above with the characterization of the pre-infimum process in Proposition \ref{prop.preinfcont}, we obtain the following result.

\begin{proposition}\label{prop.constructionLandH}
Let $\hat{L}$ be a process which under $\P$ has the law of $L$ under $\P^\#$, satisfying condition {(C3)}, so that its height process is well-defined and has a continuous modification $\hat{H}$. Define $\hat{I}_t=\inf\{\hat{L}_s:s\leq t\}$. Let $E$ be an exponential random variable with rate $\xi$. Let $\widetilde{D}$ be distributed as in the statement of Theorem \ref{thm.postinfcont}, independent of $\hat{L}$ and $E$, and set $$D_t=\widetilde{D}_t+E,$$ and
$${D}^{-1}_t:=\inf\{s:{D}_s>t\}.$$
Then the height process of $L$ is well-defined and has a continuous modification $H$. Moreover, $H$ is also the height process of $L-\underline{\underline{L}}$ and 
$$(L_t-\underline{\underline{L}}_t, H_t,t\geq 0)\overset{d}{=}(\hat{L}_t+D(D^{-1}(\hat{I}_t)), \hat{H}_t+\widetilde{D}^{-1}(\hat{I}_t), t\geq 0).$$
\end{proposition}
\begin{proof}
The existence of $\hat{H}$ follows from Proposition \ref{prop.preinfcont}. The construction of $L-\underline{\underline{L}}$ follows from Proposition \ref{prop.preinfcont} for the pre-infimum process, and from Lemma 4.1 in \cite{Millar1977} and Theorem \ref{thm.postinfcont} for the post-infimum process.

We claim that $(\hat{H}_t+\widetilde{D}^{-1}(\hat{I}_t),t\geq 0)$ is the height process corresponding to the process $$(\hat{L}_t+D(D^{-1}(\hat{I}_t)),t\geq 0).$$ Firstly, note that for $\rho^*=\sup\{t:\hat{I}_t=-E\}$,
$$(\hat{L}_t+D(D^{-1}(\hat{I}_t)) , \hat{H}_t+\widetilde{D}^{-1}(\hat{I}_t), t\in [0,\rho^*])=(\hat{L}_t+E, \hat{H}_t, t\in [0,\rho^*]).$$
The height process is, by definition, not affected by adding a constant to the \L ukasiewicz path, so on $[0,\rho^*]$, the claim follows. On $[\rho^*,\infty)$, the claim follows from Theorem \ref{thm.postinfcont}. 

Finally, we claim that the height processes of $L$ and $L-\underline{\underline{L}}$ agree. Set $I_\infty=\inf\{L_t,t\geq 0\}$, and $\rho=\sup\{t>0:L_{t}=I_\infty\}$, and observe that
$$(L_t-\underline{\underline{L}}_t , t\in [0,\rho])=(L_t-I_\infty, t\in [0,\rho]).$$
Again using the fact that the height process is not affected by adding a constant to the \L ukasiewicz path, the claim follows on $[0,\rho]$. On $[\rho,\infty)$, the claim follows from Lemma 8ii in \cite{Lambert2002}.  
\end{proof}

\section{Joint convergence of the height process and \L ukasiewicz path} \label{sec.convergence}

In this section, we will show the convergence of the discrete \L ukasiewicz path $S^n$ and height process $H^n$ to their continuous counterparts $L$ and $H$ under rescaling. The convergence result relies on the construction of the discrete and continuous processes introduced in \eqref{eq.pathwiseconstrSandHentireprocess} and Proposition \ref{prop.constructionLandH} respectively. 

We will start by proving the joint convergence of $\hat{S}^n$ and its height process $\hat{H}^n$ under rescaling. 
\begin{theorem}\label{thm.jointconvSandHpathwiseconstr}
We have that 
$$\left(n^{-1}\hat{S}^n_{\lfloor n \gamma_n t \rfloor}, \gamma_n^{-1} \hat{H}^n_{\lfloor n \gamma_n t \rfloor}, t\geq 0\right) \overset{d}{\to} \left(\hat{L}_t,\hat{H}_t,t\geq 0\right)$$
in $\D(\R_+,\R)^2$ as $n\to \infty$.
\end{theorem}
\begin{proof}
We will first show convergence under rescaling of the \L ukasiewicz path, i.e.\ 
$$\left(n^{-1}\hat{S}^n_{\lfloor n \gamma_n t \rfloor}, t\geq 0\right) \overset{d}{\to} \left(\hat{L}_t,t\geq 0\right)$$
in the Skorokhod topology as $n\to \infty$, after which we will use Theorem 2.3.1 of Duquesne and Le Gall \cite{Duquesne2005} to show the joint convergence with the height process. Since $\hat{S}^n$ is a downward skip-free random walk on the integers, it will be sufficient to show that for all $\theta>0$,
\begin{equation}\label{eq.convLukasiewiczcondtodieout}\E\left[\exp\left(-\theta n^{-1}\hat{S}^n_{\lfloor n\gamma_n t \rfloor}\right)\right]\to \E\left[\exp(-\theta \hat{L}_t)\right]=\exp(t \phi(\theta+\xi)).\end{equation}
Note that by definition,
\begin{align*}
    \E\left[\exp\left(-\theta n^{-1}\hat{S}^n_{\lfloor n \gamma_n t \rfloor}\right)\right]&=\E\left[\exp\left(-\xi_n {S}^n_{\lfloor n \gamma_n t \rfloor}\right)\exp\left(-\theta n^{-1}{S}^n_{\lfloor n\gamma_n t \rfloor}\right)\right]\\
    &=\E\left[\exp\left(-\left(n\xi_n+\theta\right)n^{-1}{S}^n_{\lfloor n \gamma_n t \rfloor}\right)\right]\\
    &=\exp\left(\lfloor n \gamma_n t \rfloor \phi_n\left(n^{-1}\left(n\xi_n+\theta\right)\right)\right)
\end{align*}
For sake of brevity, we will denote $n\gamma_n\phi_n(n^{-1}\cdot)$ by $\bar{\phi}_n(\cdot)$. By the convergence under rescaling of $S^n$ to $L$, we know that $\bar{\phi}_n$ converges to $\phi$ pointwise. We will first show that \begin{equation}\label{eq.convergencexi}n\xi_n\to\xi\end{equation} as $n\to \infty$; then we will show that there is a $b<\xi$ such that for any $B>b$,
\begin{equation}\label{eq.uniformconvphi} \bar{\phi}_n(\cdot) \to \phi(\cdot)\end{equation}
uniformly on $[b,B]$. Together, \eqref{eq.convergencexi} and \eqref{eq.uniformconvphi} imply \eqref{eq.convLukasiewiczcondtodieout}. \\
Note that by definition, $n\xi_n$ is the unique non-trivial zero of $\bar{\phi}_n$. By convexity of $\phi$ and by $\phi(0)=\phi(\xi)=0$, we see that for all $\epsilon>0$, 
$$\phi(\xi-\epsilon)<0<\phi(\xi+\epsilon).$$
So by the pointwise convergence of $\bar{\phi}_n$ to $\phi$, we have that, for all $n$ large enough, $$\bar{\phi}_n(\xi-\epsilon)<0<\bar{\phi}_n(\xi+\epsilon).$$
But then $n\xi_n\in (\xi-\epsilon,\xi+\epsilon)$,
which implies \eqref{eq.convergencexi}. To prove \eqref{eq.uniformconvphi}, given the pointwise convergence of $\bar{\phi}_n$ to $\phi$, it will be sufficient to find a $y<\xi$ such that, for all $n$ large enough, $\bar{\phi}_n$ is monotone on $[y,\infty)$. By convexity of $\bar{\phi}_n$ for all $n$, it is enough to show that there is an $x<y<\xi$ such that \begin{equation}\label{eq.alphabeta}\bar{\phi}_n(x)<\bar{\phi}_n(y)<0\end{equation} for all $n$ large enough. However, by convexity of $\phi$ and by $\phi(0)=\phi(\xi)$, there exist $x<y<0$ such that 
$$\phi(x)<\phi(y)<0.$$
The pointwise convergence of $\bar{\phi}_n$ to $\phi$ then implies \eqref{eq.alphabeta}, and hence \eqref{eq.uniformconvphi} and \eqref{eq.convLukasiewiczcondtodieout}.

We now wish to apply Theorem 2.3.1 in \cite{Duquesne2005} to obtain joint convergence under rescaling of the \L ukasiewicz path and height process. Considering Theorem \ref{thm.factsonprocessdyingout}.3, Condition {(C3)} and \eqref{eq.convLukasiewiczcondtodieout}, the only condition that is left to check is that for ${Y}^n_m$ the number of individuals in generation $m$ in the Galton--Watson branching process given by the first $n$ trees in the forest encoded by $S^n$, we have for all $\delta>0$,
$$\liminf\limits_{n\to\infty}\P^\#\left[{Y}^n_{\lfloor \delta \gamma_n\rfloor}=0\right]>0.$$
We claim that this is equivalent to the statement above under $\P$. Indeed,
$$\P\left[Y^n_{\lfloor\delta\gamma_n\rfloor}=0\right]=\P\left[
\text{the first }n\text{ trees are finite}
\right] \times \P^\#\left[{Y}^n_{\lfloor \delta \gamma_n\rfloor}=0\right],$$
so the equivalence follows from the fact that 
$$\P\left[\text{the first }n\text{ trees are finite}\right]\to \exp(-\xi).$$
The statement then follows from 
$$\liminf\limits_{n\to\infty}\P\left[{Y}^n_{\lfloor \delta \gamma_n\rfloor}=0\right]>0,$$
which is assumption \eqref{eq.assumptiontreesdie}. \end{proof}


We will now show joint convergence under rescaling of $Q^n$ and $D^n$, which is the content of the following theorem. Suppose that the Laplace exponent of $L$ is given by $$\phi(\theta)=a \theta +\frac{1}{2}b \theta^2+\int_0^\infty \nu(dx)(\exp(-\theta x)-1+\theta (x\wedge 1)).$$
Define 
$$\tilde{\nu}(dx,dy)=\exp(-\xi x)\nu(dy)dx \one_{x<y}$$ and let $(D_t,Q_t)_{t\geq 0}$ be a subordinator with L\'evy measure $\tilde{\nu}$ and drift vector $(b,2b)$. Let $E$ be an exponential random variable with rate $\xi$.

\begin{theorem}\label{thm.convergenceDnQn}
  It holds that
$$\left(n^{-1}D^n_{\lfloor \gamma_n t \rfloor}, n^{-1}Q^n_{\lfloor \gamma_n t \rfloor},t\geq 0\right)\overset{d}{\to}(E+D_t,E+Q_t, t\geq 0)$$
in the Skorokhod topology, as $n\to \infty$.
\end{theorem}
\begin{proof}
Firstly, we recall that $D^n_0=Q^n_0\sim \operatorname{Geom}(\exp(-\xi_n))$. Recalling that $n\xi_n\to \xi$, we get that for any $x>0$,
\begin{align*}\exp(-x n^{-1} D^n_0)&=\frac{1-\exp(-\xi_n)}{1-\exp(-\xi_n-x n^{-1})}\\
&\to \frac{\xi}{\xi+x}, \end{align*}
so we can conclude that 
\begin{equation}\label{eq.convergenceGtoE}(n^{-1} D^n_0, n^{-1} Q^n_0)\overset{d}{\to}(E,E).\end{equation}
We aim to use \cite[Theorem VII.2.9]{JacodShiryaev2003} to prove the convergence under rescaling of 
$$\left(\sum\limits_{i=1}^k B_i^{(n)}, \sum\limits_{i=1}^k (N_i^{(n)}-1), k\geq 1 \right).$$ We will approximate this random walk by a Poissonized version. To that end, define $$\nu_n(dx)=n\gamma_n \sum_{i=0}^\infty\delta_{(i-1)/n}(dx)\P(D^n=i).$$
We claim that for 
$(\bar{S}^n_t,t\geq 0)$ the compound Poisson process defined by 
$$\E\left[\exp(iu\bar{S}^n_t)\right]=\exp\left(t\int_{-1/n}^\infty \nu_n(dx)(e^{iux}-1)\right),$$
it holds that
\begin{equation}\label{eq.convergencepointprocesstildeS}
\bar{S}^n_1 \overset{d}{\to}L_1.
\end{equation}
Indeed, we note that $\bar{S}^n_1$ is distributed as a Poisson process with rate $1$ with jumps distributed as $n^{-1}(D^n-1)$ evaluated at time $n\gamma_n $, which is a Poissonized version of $S^n_{\lfloor n\gamma_n \rfloor}$, so \eqref{eq.convergencepointprocesstildeS} is implied by \eqref{eq.assumptionconvergence}. Now, let $h:\R\to \R$ be a bounded function such that $h(x)=x$ on a neighbourhood of $0$. Then, Theorem VII.2.9 in \cite{JacodShiryaev2003} implies the following facts. 
\begin{enumerate}
\item $\int_{-1/n}^\infty\nu_n(dx)h(x)\to a_h$ for some $a_h$ as $n\to\infty$.
\item As $n\to \infty$,
\begin{equation}\int_{-1/n}^\infty\nu_n(dx)(h(x))^2 \to b+\int_{0}^\infty\nu(dx)(h(x))^2.\end{equation}
\item  For bounded continuous $f:\R^+\to \R$ such that $f(x)=o(x^2)$ when $x\downarrow 0$, it holds that
\begin{equation}\label{eq.convergenceintegraltestfunction}\int_{-1/n}^\infty\nu_n(dx)f(x) \to \int_{0}^\infty\nu(dx)f(x).\end{equation}
\end{enumerate}
Recalling that $$\P(B^{(n)}=k, N^{(n)}=l)=\exp(-k\xi_n)\P(D=l)\one_{k<l},$$ 
we now define 
$$\tilde{\nu}_n(dx,dy)=\gamma_n\sum_{j=0}^\infty\sum_{i=0}^{j-1}\delta_{i/n,j/n}(dx,dy)\exp(-\xi_n i)\P(D^n=j).$$
Then, by a similar argument to before, for 
$(\bar{D}^n_t, \bar{Q}^n_t,t\geq 0)$ the process with stationary increments defined by 
$$\E\left[\exp(iu_1\bar{D}^n_t+iu_2\bar{Q}^n_t)\right]=\exp\left(t\int_{0}^\infty\int_{0}^\infty \tilde{\nu}_n(dx,dy)(e^{iu_1x+iu_2 y}-1)\right),$$
we get that for all $t$, $(\bar{D}^n_t, \bar{Q}^n_t)$ is a Poissonised version of $$\left(n^{-1}\sum_{i=1}^{\lfloor\gamma_nt\rfloor} B_i^n,n^{-1}\sum_{i=1}^{\lfloor\gamma_nt\rfloor} N_i^n\right).$$
We will show that for all $t$,
\begin{equation}\label{eq.convDandQ} (\bar{D}^n_t, \bar{Q}^n_t)\overset{d}{\to}(D_t,Q_t),\end{equation}
which implies the functional convergence by a result in Kallenberg's book \cite[Theorem 16.14]{Kallenberg2002}.
By Theorem VII.2.9 in \cite{JacodShiryaev2003}, using truncation function $H(x,y)=(x\wedge 1,y\wedge 1)$, the following properties imply \eqref{eq.convDandQ}.
\begin{enumerate}
    \item \begin{enumerate}
        \item \label{enumconvintegrals1a}$\int_0^\infty\int_0^\infty \tilde{\nu}_n(dx,dy)(x\wedge 1)\to \int_0^\infty\int_0^\infty \tilde{\nu}(dx,dy)(x\wedge 1)+b$
        \item \label{enumconvintegrals1b}$\int_0^\infty\int_0^\infty \tilde{\nu}_n(dx,dy)(y\wedge 1)\to \int_0^\infty\int_0^\infty \tilde{\nu}(dx,dy)(y\wedge 1)+2b$
        \end{enumerate}
    \item \begin{enumerate}
        \item \label{enumconvintegrals2a}$\int_0^\infty\int_0^\infty \tilde{\nu}_n(dx,dy)(x\wedge 1)^2\to \int_0^\infty\int_0^\infty \tilde{\nu}(dx,dy)(x\wedge 1)^2$
        \item \label{enumconvintegrals2b}$\int_0^\infty\int_0^\infty \tilde{\nu}_n(dx,dy)(y\wedge 1)^2\to \int_0^\infty\int_0^\infty \tilde{\nu}(dx,dy)(y\wedge 1)^2$
        \item \label{enumconvintegrals2c}$\int_0^\infty\int_0^\infty \tilde{\nu}_n(dx,dy)(x\wedge 1)(y\wedge 1)\to \int_0^\infty\int_0^\infty \tilde{\nu}(dx,dy)(y\wedge 1)^2$ \end{enumerate}
    \item \label{enumconvintegrals3}For all continuous, bounded $F:\R^2\to \R$ that are $0$ on a neighbourhood of $0$, 
    $$\int_0^\infty\int_0^\infty \tilde{\nu}_n(dx,dy)F(x,y)\to \int_0^\infty\int_0^\infty \tilde{\nu}(dx,dy)F(x,y).$$
\end{enumerate}
We will prove the conditions one-by-one, starting with \ref{enumconvintegrals1a}. We note that
\begin{align*}\int_0^\infty\int_0^\infty \tilde{\nu}_n(dx,dy)(x\wedge 1)&=\gamma_n\sum_{i=0}^\infty\sum_{j=0}^{i-1}\P(D^n=i)\exp(-\xi_n j)(j/n\wedge 1)\\
&= \frac{1}{n}\int_{-1/n}^\infty\nu_n(dx) \sum_{j=0}^{\lfloor nx \rfloor }\exp(-\xi_n j)(j/n\wedge 1)dy. \end{align*}
We will first argue that 
\begin{equation}\label{eq.uniformconvergenceriemannint}\frac{1}{n(x\wedge 1)} \sum_{j=0}^{\lfloor nx\rfloor }\exp(-\xi_n j)(j/n\wedge 1)\to \frac{1}{1\wedge x}\int_0^x dy\exp(-\xi y)(y\wedge 1)\end{equation}
uniformly in all $x>0$ as $n\to \infty$. Recall that $n\xi_n=\xi+o(1)$, which, together with the Riemann integrability of $\exp(-\xi y)(y\wedge 1)$ on compact intervals implies pointwise convergence in \eqref{eq.uniformconvergenceriemannint}. The facts that $\int_0^\infty dy\exp(-\xi y)(y\wedge 1)<\infty$ and that $\int_0^x dy\exp(-\xi y)(y\wedge 1)=O(x^2)$ as $x\to 0$, together with the uniform continuity of $f(x,y):=\frac{1}{x\wedge 1}\exp(-\xi y)(y\wedge 1)$ on $(\delta,\infty)\times (0,R)$ for any $\delta>0$ and $R>0$, then imply that convergence in \eqref{eq.uniformconvergenceriemannint} is in fact uniform, as claimed. \\
Since $\int_{-1/n}^\infty\nu_n(dx)(x\wedge 1)$ is convergent as $n\to \infty$, and in particular bounded, we have that
\begin{equation}\label{eq.convergenceriemannintegrated}\left|\int_{-1/n}^\infty\nu_n(dx) \frac{1}{n}\sum_{j=0}^{\lfloor nx \rfloor}\exp(-\xi_n j)(j/n\wedge 1)dy-\int_{-1/n}^\infty\nu_n(dx)\int_0^x dy\exp(-\xi y)(y\wedge 1)\right|\to 0\end{equation}  as $n\to \infty$. So in order to prove \ref{enumconvintegrals1a}, it is sufficient to show that 
$$\int_{-1/n}^\infty\nu_n(dx)\int_0^x dy\exp(-\xi y)(y\wedge 1)\to \int_0^\infty\int_0^\infty \tilde{\nu}(dx,dy)(x\wedge 1)+b.$$
 To see this, consider a truncation function $h$, and define $$g_u(x):=e^{iux}-1+iuh(x)-u^2\int_0^x dy\exp(-\xi y)(y\wedge 1).$$ 
 Then, note that $\int_0^x dy\exp(-\xi y)(y\wedge 1)=\frac{x^2}{2}+o(x^2)$ as $x\to 0$, so $g_u(x)=o(x^2)$ as $x\to 0$. Furthermore, $g_u$ is bounded, because $\int_0^x dy\exp(-\xi y)(y\wedge 1)$ is bounded. Then, by \eqref{eq.convergenceintegraltestfunction}, 
\begin{equation}\label{eq.convergenceg_u(x)}\int_{-1/n}^\infty \nu_n(dx)g_u(x)\to \int_0^\infty \nu(dx) g_u(x).\end{equation}
Moreover, since 
\begin{align*}&\E\left[\exp(iu\bar{S}^n)\right]=\exp\left(-iu\int_{-1/n}^\infty \nu_n(dx)h(x)+\int_{-1/n}^\infty \nu_n(dx)(e^{iux}-1+uih(x))\right)\\
&\to \E\left[\exp(iuL_1)\right]
=\exp\left(-iua_h-b u^2+\int_{0}^\infty \nu(dx)(e^{iux}-1+uih(x))\right)\end{align*}
and $\int_{-1/n}^\infty \nu_n(dx)h(x)\to a_h$, we find that $$\int_{-1/n}^\infty \nu_n(dx)(e^{iux}-1+uih(x))\to {b}u^2+\int_{0}^\infty \nu(dx)(e^{iux}-1+uih(x)),$$
which, combined with \eqref{eq.convergenceg_u(x)}, implies that 
$$\int_{-1/n}^\infty\nu_n(dx)\int_0^x dy\exp(-\xi y)(y\wedge 1)\to \int_0^\infty\int_0^\infty \tilde{\nu}(dx,dy)(x\wedge 1)+b,$$ proving \ref{enumconvintegrals1a}.\\
A similar proof shows that also \ref{enumconvintegrals1b} holds.\\
To prove \ref{enumconvintegrals2a}, note that 
\begin{align*}\int_0^\infty\int_0^\infty \tilde{\nu}_n(dx,dy)(x\wedge 1)^2&=\gamma_n\sum_{i=0}^\infty\sum_{j=0}^{i-1}\P(D^n=i)\exp(-\xi_n j)(j/n\wedge 1)^2\\
&= \frac{1}{n}\int_{-1/n}^\infty\nu_n(dx) \sum_{j=0}^{\lfloor nx \rfloor}\exp(-\xi_n j)(j/n\wedge 1)^2dy. \end{align*}
The first fact we need is that 
$$\left|\frac{1}{n}\int_{-1/n}^\infty\nu_n(dx) \sum_{j=0}^{\lfloor nx \rfloor}\exp(-\xi_n j)(j/n\wedge 1)^2dy-\int_{-1/n}^\infty\nu_n(dx)\int_0^x dy\exp(-\xi y)(y\wedge 1)^2\right|\to 0$$ as $n\to \infty$,
which is proved in a similar manner to \eqref{eq.convergenceriemannintegrated}. Then, since 
$\int_0^x dy\exp(-\xi y)(y\wedge 1)^2=o(x^2)$ as $x\to 0$, and it is a bounded function of $x$, we can use \eqref{eq.convergenceintegraltestfunction} to obtain \ref{enumconvintegrals2a}. The proofs of \ref{enumconvintegrals2b}, \ref{enumconvintegrals2c}, and \ref{enumconvintegrals3} are similar. \\
By Theorem VII.2.9 in \cite{JacodShiryaev2003}, for all $t>0$,
\begin{equation*} (\bar{D}^n_t, \bar{Q}^n_t)\overset{d}{\to}(D_t,Q_t),\end{equation*}
which proves the statement. \end{proof}
The final steps in the proof of Theorem \ref{thm.mainconvergenceresult} are similar to the proof of Theorem 1.5 in \cite{Duquesne2009}. 

Recall that $F^{n}(k)=\inf\{l:-\hat{I}^n(k-1-l)<D^n(l)\}$. We will now use the convergence of $(D^n,Q^n)$ under rescaling to prove the convergence under rescaling of $F^{n}$. For this, we need the following technical lemma, whose proof may be found in the appendix.

\begin{lemma}\label{lemma.technicalconvergenceinfimum}
Suppose $f_n\to f$ in $\D(\R_+,\R)$ as $n\to\infty$ with $f_n$ increasing for all $n$, and $f$ increasing and continuous. Furthermore, suppose $g_n\to g$ in $\D(\R_+,\R)$ as $n\to\infty$ with $g_n$ increasing for all $n$, and $g$ strictly increasing, and $g_n(t) \to \infty$ and $g(t)\to \infty$ as $t\to \infty$. Then,
$$\left(\inf \{s\geq 0:f(t)<g(s)\},t \geq 0 \right)$$
is continuous, and
$$\left(\inf \{s\geq 0:f_n(t)<g_n(s)\},t \geq 0 \right)\to \left(\inf \{s\geq 0:f(t)<g(s)\},t \geq 0 \right)$$
in $\D(\R_+,\R)$  as $n\to\infty$.\\
Moreover, for $(\epsilon_n)_{n\geq 0}$ such that $\epsilon_n\downarrow 0$,
$$\left(\inf \{s\geq 0:f_n(t-\epsilon_n s)<g_n(s)\},t \geq 0 \right)\to \left(\inf \{s\geq 0:f(t)<g(s)\},t \geq 0 \right)$$
in $\D(\R_+,\R)$ as $n\to\infty$.
\end{lemma}

\begin{lemma}\label{lemma.convFn}
We have that
$$\left(\gamma_n^{-1}F^{n}\left(\lfloor n \gamma_n  t \rfloor\right),t\geq 0\right) \overset{d}{\to} \left(F_t,t\geq 0\right)$$
in $\D(\R_+,\R)$ as $n\to\infty$, jointly with the joint convergence of $D^n$ and $Q^n$ under rescaling (Theorem \ref{thm.convergenceDnQn}) and the joint convergence of $\hat{S}^n$ and $\hat{H}^n$ under rescaling (Theorem \ref{thm.jointconvSandHpathwiseconstr}), where $(D^n,Q^n)$ is independent of $(\hat{S}^n,\hat{H}^n)$. Moreover, $(F_t,t\geq 0)$ is continuous almost surely.
\end{lemma}
\begin{proof}
We need to show that
$$\left(\gamma_n^{-1}\inf\{l:-\hat{I}^n(\lfloor n \gamma_n t \rfloor-1-l)<D_l\},t\geq 0\right) \overset{d}{\to} \left(\inf\{s:-\hat{I}_t<D_s\},t\geq 0\right)$$
in $\D(\R_+,\R)$ as $n\to\infty$ and that $\left(\inf\{s:-\hat{I}_t<D_s\},t\geq 0\right)$ is continuous, for which we will use Lemma \ref{lemma.technicalconvergenceinfimum}. Firstly, recall that $L$ is of infinite variation, so that $b>0$, or $\nu((0,\infty))=\infty$. In the first case, it is obvious that $D$ is strictly increasing. In the second case, note that since $D$ has L\'evy measure 
$$\int_{y\in[0,\infty)} \tilde{\nu}(dx,dy)=\exp(-\xi x)\nu((x,\infty))dx,$$
the intensity of jumps of size $\epsilon$ goes to $\infty$ as $\epsilon$ goes to zero, which implies that the jumps of $D$ are dense, and its sample paths are strictly increasing with probability $1$.
By Skorokhod's representation theorem, we may work on a probability space where the joint convergence of $D^n$ and $Q^n$ under rescaling (Theorem \ref{thm.convergenceDnQn}) and the joint convergence of $\hat{S}^n$ (and $\hat{I}^n$) and $\hat{H}^n$ under rescaling (Theorem \ref{thm.jointconvSandHpathwiseconstr}) all hold almost surely. Then, by Lemma \ref{lemma.technicalconvergenceinfimum}, since $\hat{I}$ is continuous ($\hat{L}$ is spectrally positive), 
$$\left(\gamma_n^{-1}\inf\{l:-\hat{I}^n(\lfloor \gamma_nnt\rfloor-1-l)<D_l\},t\geq 0\right) \overset{\text{a.s.}}{\to} \left(\inf\{s:-\hat{I}_t<D_s\},t\geq 0\right)$$
in $\D(\R_+,\R)$ as $n\to\infty$ and $F_t$ is continuous almost surely. The result follows. 
\end{proof}

From Theorem \ref{thm.jointconvSandHpathwiseconstr}, Theorem \ref{thm.convergenceDnQn},  and Lemma \ref{lemma.convFn} we know that, as $n\to \infty$,
\begin{align}\begin{split}\label{eq.jointconvergenceeverything}
   \left( n^{-1} \hat{S}^{n}\left(\lfloor n \gamma_n t\rfloor\right), \gamma_n^{-1}\hat{H}^{n}\left(\lfloor n \gamma_n t\rfloor\right), t\geq 0\right)&\overset{d}{\to} \left(\hat{L}_t,\hat{H}_t, t\geq 0\right),\\
   \left(n^{-1}D^n\left(\lfloor \gamma_n t\rfloor\right), n^{-1}Q^n\left(\lfloor \gamma_n t\rfloor\right), t\geq 0\right)&\overset{d}{\to}\left(D_t,Q_t,t\geq 0\right)\text{, and}\\
   \left(\gamma_n^{-1}F^{n}\left(\lfloor n \gamma_n t \rfloor\right), t\geq 0\right)&\overset{d}{\to}\left(F_t,t\geq 0\right),
\end{split}
\end{align}
jointly, in $\D(\R_+,\R)^2$, $\D(\R_+,\R)^2$, and $\D(\R_+,\R)$ respectively.
We would like to show the convergence under rescaling of $Q^n(F^{n}(k))$ and of $D^n(F^{n}(k))$ jointly with the convergence in \eqref{eq.jointconvergenceeverything}, since these quantities appear in the pathwise construction of $S^n$ and $S^n-\underline{\underline{S}}^n$ in equation \eqref{eq.pathwiseconstrSandHentireprocess}. For this we need a technical lemma, which follows immediately from the characterization of convergence in the Skorokhod topology given in the book of Ethier and Kurtz \cite[Proposition 3.6.5]{Ethier1986}. 
\begin{lemma}\label{lemma.technicalcomposedfunctions}
If $h_n\to h$ and $f_n\to f$ in $\D(\R_+,\R)$ as $n\to\infty$, $h_n$ and $h$ are monotone non-decreasing, and $h$ is continuous, then 
$$h_n\circ f_n \to h\circ f$$
in $\D(\R_+,\R)$ as $n\to\infty$.
\end{lemma}

\begin{lemma}\label{lemma.convsecondterms}
We have 
$$\left(n^{-1}D^n\left(F^{n}\left(\lfloor n \gamma_n t\rfloor\right)\right),n^{-1}Q^n\left(F^{n}\left((\lfloor n \gamma_n t\rfloor\right)\right), t\geq 0 \right)\overset{d}{\to} \left(D_{F_t}, Q_{F_t},t\geq 0\right),$$
in $\D(\R_+,\R^2)$ as $n\to\infty$, jointly with the convergence in \eqref{eq.jointconvergenceeverything}. 
\end{lemma}
\begin{proof}
By Skorokhod's representation theorem we may work on a space where the convergence of \eqref{eq.jointconvergenceeverything} holds almost surely. The result then follows from Lemma \ref{lemma.technicalcomposedfunctions}. 
\end{proof}
\begin{lemma}\label{lemma.convfirstterms}
We have 
\begin{align*}
&\left(n^{-1}\hat{S}^n\left(\lfloor n \gamma_n t\rfloor -F^{n}\left(\lfloor n \gamma_n t\rfloor\right)\right), \gamma_n^{-1}\hat{H}^n\left(\lfloor n \gamma_n t\rfloor -F^{n}\left(\lfloor n \gamma_n t\rfloor\right)\right),t\geq 0 \right)\\&\overset{d}{\to} (\hat{L}_t,\hat{H}_t \geq 0)\end{align*}
in $\D(\R_+,\R)^2$ as $n\to\infty$ jointly with the convergence in \eqref{eq.jointconvergenceeverything}. 

\end{lemma}
\begin{proof}
Firstly, note that $k-F^{n}(k)$ equals the number of steps not spent on the spine up to time $k$ and so is a non-decreasing function of $k$. Then, note that
\begin{align*}
    \lfloor n \gamma_nt\rfloor -F^{n}\left(\lfloor n \gamma_nt\rfloor\right) &=  \left\lfloor n \gamma_n\left(t -n^{-1} F^{n}\left(\lfloor n \gamma_nt\rfloor\right)\right)\right\rfloor,
\end{align*}
and 
$$\left(t-\gamma_n^{-1}n^{-1}F^{n}\left(\lfloor n \gamma_nt\rfloor\right),t\geq 0\right )\to (t,t\geq 0)$$
in $\D(\R_+,\R)$ almost surely as $n\to\infty$. 
We may use Skorokhod's representation theorem to work on a space where the convergence in \eqref{eq.jointconvergenceeverything} holds almost surely, and then Lemma \ref{lemma.technicalcomposedfunctions} gives the result.  \end{proof}

We will now prove Theorem \ref{thm.mainconvergenceresult}. 

\begin{proof}[Proof of Theorem \ref{thm.mainconvergenceresult}]
Let $(\hat{L}, \hat{H})$, $\hat{I}$, $(D,Q)$, and  $F_t=\inf\{s:-\hat{I}_t<D_s\}$ be as defined in Section \ref{sec.supercriticalheightprocess}.
Then, for $\underline{\underline{S}}^n$ the future infimum of $S^n$, we know by the pathwise construction of $S^n$, $S^n-\underline{\underline{S}}^n$ and $H^n$  given in \eqref{eq.pathwiseconstrSandHentireprocess}, Lemmas \ref{lemma.convfirstterms}, \ref{lemma.convsecondterms} and \ref{lemma.convFn}, that 
\begin{align}\begin{split}\label{eq.pathwiseconvergences} &\left(n^{-1} S^n\left(\lfloor tn \gamma_n\rfloor\right),n^{-1} (S^n-\underline{\underline{S}}^n)\left(\lfloor tn \gamma_n\rfloor\right), \gamma_n^{-1} H^{n}\left(\lfloor tn \gamma_n\rfloor\right), t\geq 0\right)\\
&\xrightarrow{d} \left(\hat{L}_t+Q_{F_t},\hat{L}_t+D_{F_t},\hat{H}_t+F_t, t\geq 0\right)\end{split}\end{align}
in $\D(\R_+,\R^3)$ as $n\to\infty$ .
By assumption, we have that 
$$\left(  n^{-1}S^n\left(\lfloor n \gamma_n t\rfloor\right), t\geq 0\right)\overset{d}{\to}(L_t,t\geq 0),$$
so 
$$L_t\overset{d}{=}\hat{L}_t+Q_{F_t},$$ and by construction, $\hat{L}_t+D_{F_t}$ equals $\hat{L}_t+Q_{F_t}$ minus its future infimum. Then, by Proposition \ref{prop.constructionLandH}, we know that $\hat{H}_t+F_t$ is the height process corresponding to $\hat{L}_t+D_{F_t}$ and hence to $\hat{L}_t+Q_{F_t}.$ The result follows.
\end{proof}

\section{Application to the configuration model in the critical window}\label{sec.application}

This section contains new results on the scaling limit of the configuration model with i.i.d.\ power-law degrees in the critical window. We use Theorem \ref{thm.mainconvergenceresult} to extend the methods in Conchon-Kerjan and Goldschmidt \cite{Conchon2018} from the critical point  to the critical window.

The configuration model is a method to construct a multigraph with a given degree sequence that was introduced by Bollobás in \cite{Bollobas1980configurationmodel}.
\begin{quote}
Consider $n$ vertices labelled by $[n]$ and a sequence $\mathbf{d}=(d_i)_{i\in[n]}\in \N^n$ such that $\sum_{i\in [n]}d_i$ is even. We will sample a multigraph such that the degree of vertex $i$ is equal to $d_i$ for every $i\in [n]$. The configuration model on $n$ vertices having degree sequence $\mathbf{d}$ is constructed as follows. Equip vertex $j$ with $d_j$ half-edges. Two half-edges create an edge once they are paired. Pick any half-edge and pair it with a uniformly chosen half-edge from the remaining unpaired half-edges and keep repeating the above procedure until all half-edges are paired.\end{quote}
Note that the graph constructed by the above procedure may contain self-loops or multiple edges. It can be shown that, conditionally on the constructed multigraph being simple, the law of such graphs is uniform over all possible simple graphs with degree sequence $\mathbf{d}$. Furthermore, as shown in  \cite{Janson2009TPT}, under very general assumptions, the asymptotic probability of the graph being simple is positive. For a discussion of the configuration model and standard results, see \cite[Chapter 7]{VanderHofstad2017}. 

\subsection{Model and result}\label{subsec.model}
We use the configuration model to construct a uniform graph with a random degree sequence. The model we consider is as follows. 

Fix $\lambda\in\R$. Most quantities that will be defined depend on $\lambda$. To avoid overcomplicating the notation, this will not be made explicit unless necessary to avoid confusion.
For each $n \in \N$, let $D_1^{n}, D_2^{n}, \dots, D_n^{n} \geq 1 $ be an i.i.d.\ degree sequence satisfying the following properties, labeled by `CM' for `configuration model'. 
\begin{description} 
\item[(CM1)]For $\mu_n := \E[D_1^{n}] $, we have $\mu_n\rightarrow \mu$ as $n \rightarrow \infty$, with $\mu$ not depending on $\lambda$;
\item[(CM2)] $\E \left[\left(D_1^{n}\right)^2\right]=\left(2+\lambda n ^{-{(\alpha-1)/(\alpha+1)}}\right) \E[D_1^{n}]$ for some $\alpha \in (1,2)$;
\item[(CM3)]$\P (D_1^{n}=k)\sim c_n k^{-(\alpha+2)}$ as $k\rightarrow \infty$, with $c_n>0$ for all $n$ and $c_n\to c$ as $n\to \infty$;
\item[(CM4)] For $Z^n$ a random variable such that $\P(Z^n=k)=k\P(D_1^{n}=k)/\mu_n$, and $(S^n(k),k\geq 0)$ a random walk with steps distributed as $Z^n-2$, let $Y^n_m$ be the number of vertices at height $m$ in the first $\lfloor n^{1/(\alpha+1)}\rfloor$ trees of the forest encoded by $(S^n(k),k\geq 0)$. Then, for every $\delta>0$, \begin{equation*}\liminf\limits_{n\to\infty}\P\left[Y^n_{\lfloor \delta n^{(\alpha-1)/(\alpha+1)}\rfloor}=0\right]>0.\end{equation*}
\end{description}
Let $G_1^n, G_2^n,\dots$ be the components of a uniformly random graph with i.i.d.\ degrees that are distributed as $D_1^n$, with the components listed in decreasing order of size. View each $G_i^n$ as a compact measured metric spaces by equipping it with the graph distance $d_i^n$, and the counting measure on its vertices, $\mu_i^n$. More generally, compact measured metric spaces will be denoted by a triple $(G,d,\mu)$, for $(G,d)$ a compact metric space and $\mu$ a finite Borel measure on $(G,d)$. Formally, each $G_i^n$ is an element of the Polish space of isometry-equivalence classes of measured metric spaces, endowed with the Gromov--Hausdorff--Prokhorov distance. For a discussion of the topology, we refer the reader to \cite[Section 2]{AddarioBerry2017}. We will prove the following theorem.
\begin{theorem}\label{thm.convgraphscriticalwindow}
There exists a sequence of random compact measured metric spaces $$((\mathcal{G}_1,d_1,\mu_1),(\mathcal{G}_2,d_2,\mu_2),\dots)$$ such that, as $n\to \infty$,
$$\left(\left(G_i^n,n^{-(\alpha-1)/(\alpha+1)}d_i^n, n^{-\alpha/(\alpha-1)}\mu_i^n\right),i\geq 1\right)\overset{d}{\to}((\mathcal{G}_i,d_i,\mu_i),i\geq 1)$$
in the sense of the product Gromov--Hausdorff--Prokhorov topology. \end{theorem}
If $\lambda=0$, and the degree distribution does not depend on $n$, this was already known from  \cite[Theorem 1.1]{Conchon2018}; this is known as the critical case. Intuitively, criticality entails that for large $n$, and for $(V_n,W_n)$ an edge chosen uniformly at random from the graph, the expected degree of $V_n$ is roughly $2$. Our contribution is then to prove the theorem for all $\lambda\in\R$ and for degree distributions depending on $n$. This is known as the \emph{critical window}, in which, for large $n$, and for $(V_n,W_n)$ an edge chosen uniformly at random from the graph, the expected degree of $V_n$ is roughly $2+ \lambda n ^{-{(\alpha-1)/(\alpha+1)}}$. 

\subsection{Related work}
Most results on the configuration model are obtained for models with a deterministic degree sequence. The phase transition for the undirected setting was shown in \cite{Molloy1995, Molloy1998, Janson2009}. The limiting law of the rescaled component sizes at criticality and in the critical window were obtained by Riordan \cite{Riordan2012} under the assumption that the degrees are bounded. Dhara, van der Hofstad, van Leeuwaarden and Sen showed convergence under rescaling of the component sizes and surpluses in the critical window in the finite third moment setting \cite{Dhara2017} and in the heavy-tailed regime \cite{Dhara2020}.  Bhamidi, Dhara, van der Hofstad and Sen obtained metric space convergence in the critical window in \cite{Bhamidi2020}, a result that the authors later improved to a stronger topology in \cite{Bhamidi2020Glmb}; both of these results also hold conditional on the constructed multigraph being simple. We will discuss the results by Bhamidi, Dhara, van der Hofstad and Sen further at the end of this subsection. 

Configuration models with a \emph{random} degree sequence are considered in \cite{Joseph2014} and \cite{Conchon2018}. Joseph \cite{Joseph2014} showed convergence of the component sizes and surpluses of the large components under rescaling at criticality, both for degree distributions with finite third moments and for the heavy-tailed regime. Conchon-Kerjan and Goldschmidt \cite{Conchon2018} show product Gromov--Hausdorff--Prokhorov convergence of the large components at criticality in these two regimes, and show that results also hold conditioned on the resulting graph being simple. 

In recent work \cite{Donderwinkel2021digraphs}, the author and Xie show metric space convergence under rescaling of the strongly connected components of the directed configuration model at criticality. 

\subsubsection{Results in \cite{Bhamidi2020,Bhamidi2020Glmb}}
We now describe the model that is considered in \cite{Bhamidi2020,Bhamidi2020Glmb} to provide a comparison with our result. Our description of the results in \cite{Bhamidi2020} closely follows the presentation in \cite{Conchon2018}.

 Let $\mathfrak{d}_1^n\geq \cdots\geq \mathfrak{d}_n^n$ be a family of deterministic degree sequences, such that $\sum_{i=1}^n\mathfrak{d}_i^n$ is even, and if $\mathfrak{D}_n$ denotes the degree of a vertex chosen uniformly at random, the following conditions hold, labeled by `DD' for `deterministic degrees'.
\begin{description}
    \item[(DD1)] $n^{-1/(\alpha+1)}\mathfrak{d}_i^n\to \mathfrak{\vartheta}_i$ as $n\to\infty$ for each $i\geq 1$, where $\vartheta_1\geq \vartheta_2\geq \cdots \geq 0 $ is such that $\sum_{i\geq 1} \vartheta^3<\infty$, but $\sum_{i\geq 1}\vartheta_i^2=\infty$;
    \item[(DD2)] $\mathfrak{D}_n\overset{d}{\to}\mathfrak{D}$, along with the convergence of its first two moments, for some random variable $\mathfrak{D}$ with $\P(\mathfrak{D}=1)>0)$, $\E[\mathfrak{D}]=\mu$ and $\E[\mathfrak{D}(\mathfrak{D}-1)]/\E[\mathfrak{D}]=\theta>1$, and 
    $$\lim_{K\to\infty}\limsup_{n\to\infty}n^{-3/(\alpha+1)} \sum_{i\geq K+1}\left(\mathfrak{d}_i^n\right)^3=0.$$
\end{description}
Let $\theta_n=\E\left[\mathfrak{D}_n(\mathfrak{D}_n-1)\right]/\E\left[\mathfrak{D}_n\right]$. Then, the authors sample a uniform graph with this degree sequence, and perform percolation at parameter $$p_n(\lambda)=\frac{1}{\theta_n}+\tilde{\lambda} n^{-(\alpha-1)/(\alpha+1)},$$
for some $\tilde{\lambda}\in \R$, which yields a graph in the critical window. Call the resulting degree sequence $(\mathfrak{D}_1^{n,\tilde{\lambda}}, \dots ,\mathfrak{D}_n^{n,\tilde{\lambda}})$. In this setting, \cite[Theorem 2.2]{Bhamidi2020} is the precise analogue of our Theorem  \ref{thm.convgraphscriticalwindow} in the Gromov-weak topology. In \cite{Bhamidi2020Glmb}, this result is strengthened to convergence in the Gromov--Hausdorff--Prokhorov topology, under the following additional assumptions.
\begin{description}
    \item[(DD3)] For $\mathfrak{D}^*_n$ the degree of a size-biased pick from $\mathfrak{d}_1^n,\dots,\mathfrak{d}_n^n$, there exists $c_0>0$ such that for all $1\leq l\leq \mathfrak{d}^n_1$ and $n\geq 1$,
    $$\P\left(\mathfrak{D}^*_n\geq l\right)\geq \frac{c_0}{l^\alpha}.$$
    \item[(DD4)] For $\beta_i^n=n^{-2/(\alpha+1)}\sum_{j=1}^{i-1}\left(\mathfrak{d}^n_j\right)^2$, there exists a sequence $(k_n)_{n\geq 1}$ with $k_n\to \infty$, and $k_n=o\left(n^{1/(\alpha+1)}\right)$ such that 
    $\beta_{k_n}^n=\omega(\log n)$. 
    \item[(DD5)] For all $\epsilon>0$, 
    $$\lim_{k\to \infty}\limsup_{n\to\infty} n^{-1/(\alpha+1)}\sum_{i\in (k,kn)}e^{-\epsilon \beta_i^n}\mathfrak{d}^n_i=0.$$
\end{description}
These extra assumptions allow the authors of \cite{Bhamidi2020Glmb} to show that the components in their graph model satisfy the \emph{global lower mass-bound property} \cite[Theorems 1.2 and 1.3]{Bhamidi2020Glmb}, which allows them to extend the results in \cite{Bhamidi2020} to convergence in the Gromov--Hausdorff--Prokhorov topology using the results from \cite{Athreya2016}. 

The limit object in \cite{Bhamidi2020,Bhamidi2020Glmb} is constructed by making vertex identifications in \emph{tilted inhomogeneous continuum random trees}. The scaling limit of the depth-first-walk that Bhamidi, Dhara, van der Hofstad and Sen consider is a \emph{thinned Lévy process}, whereas we show convergence to a measure changed Lévy process. The connection to their results will become clear in the following subsection.


\subsubsection{Relation to percolation}\label{subsec.percolation}
We will illustrate that the law of a degree sequence that is obtained by bond percolation on the \emph{half-edges} of a sequence of vertices with i.i.d.\ degrees in the supercritical regime satisfies the conditions of Theorem \ref{thm.convgraphscriticalwindow}. This is approximately the degree distribution after bond percolation on the \emph{edges} of a uniform random graph with such a supercritical random degree sequence, although we ignore dependence between the degrees of different vertices. Using results of Janson \cite{Janson2009a}, such mild dependence can be shown to have a negligible effect on the properties of the graph, but we omit the straightforward details. We also show how our results are related to the results in \cite{Bhamidi2020,Bhamidi2020Glmb}.


Let $D$ be a random variable in $\N$ that satisfies $\E [D]=\mu$, $\E [D^2]=\rho>2\mu$ and $\P (D=k)\sim ck^{-(\alpha+2)}$ for $\alpha \in (1,2)$ as $k\rightarrow \infty$. Define $C_{\alpha}=\frac{c\Gamma(2-\alpha)}{\alpha(\alpha-1)}$, so that the Laplace transform $\cL_D$ of $D$ satisfies
$$\cL_D(\theta)=\E [\exp (-\theta D)]=1-\mu \theta +\frac{\rho}{2}\theta^2-\frac{C_{\alpha}}{\alpha+1} \theta^{\alpha+1}+o(\theta^{\alpha+1})$$ for $\theta \rightarrow 0$. View $D$ as the degree distribution. Then, keep every half-edge with probability 
$$p(\lambda,n)=\frac{1+\lambda n^{-{(\alpha-1)/(\alpha+1)}}}{\frac{\rho}{\mu}-1},$$ and call the resulting degree distribution $B^{(\lambda,n)}$. 

In the next paragraph, we will show that the conditions of Theorem \ref{thm.convgraphscriticalwindow} are satisfied for $B^{(\lambda,n)}$. Then, for $\mathfrak{d}_1^n,\dots,\mathfrak{d}_n^n$ a sample of i.i.d.\ random variables $\mathfrak{D}_1,\dots, \mathfrak{D}_n$ with the same distribution as $D$, conditions (DD1) and (DD2) are satisfied almost surely for some sequence of random variables  $\mathfrak{\vartheta}_i$ \cite[Section 2.2]{Dhara2020}. Moreover, the order statistics of the percolated degree sequence with $\tilde{\lambda}=\lambda/(\rho/\mu-1)$ closely resemble an ordered sample of i.i.d.\ random variables distributed as $B^{(\lambda,n)}$. Therefore, it should be the case that, in this particular set-up, Bhamidi, Dhara, van der Hofstad and Sen's limit corresponds to the limit in Theorem \ref{thm.convgraphscriticalwindow}.

We will now verify the conditions of Theorem \ref{thm.convgraphscriticalwindow} for $B^{(\lambda,n)}$. Note that the Laplace transform $\cL_{B^{(\lambda,n)}}$ of $B^{(\lambda,n)}$ satisfies 
\begin{equation}\label{eq.laplacebinomial}
        \cL_{B^{(\lambda,n)}}(\theta)=\E [\exp(D\log ( 1-p(\lambda,n)+p(\lambda,n)e^{-\theta}))],
\end{equation}
such that
 \begin{align*}
     \cL_{B^{(\lambda,n)}}(\theta)
    &=1-p(\lambda,n)\mu\theta+p(\lambda,n)\mu\left(2+\lambda n^{-{(\alpha-1)/(\alpha+1)}}\right)\frac{\theta^2}{2}-\frac{C_{\alpha}}{\alpha+1}p(\lambda,n)^{\alpha+1}\theta^{\alpha+1}+o(\theta^{\alpha+1}) 
\end{align*}
as $\theta\to 0$. This implies that conditions ({CM1}), ({CM2}) and ({CM3}) are satisfied with $\mu_n=p(\lambda,n)\mu\to\mu/(\frac{\rho}{\mu}-1)$ and  $c_n=cp(\lambda,n)^{\alpha+1}\to c/(\frac{\rho}{\mu}-1)^\alpha$ as $n\to \infty$. 

We now check condition ({CM4}). We will drop the dependency on $\lambda$ from the notation, unless it is necessary to avoid confusion. Let $\tilde{D}$ be a random variable with the size-biased distribution of $D$, i.e.\ for all $k\in \N$, $$\P(\tilde{D}=k)=\frac{k\P(D=k)}{\E[D]},$$ and similarly, let $\tilde{B}_n$ be a random variable with the size-biased distribution of ${B}^{(n,\lambda)}$, i.e.\  $$\P(\tilde{B}_n=k)=\frac{k\P({B}^{(n,\lambda)}=k)}{\E\left[{B}^{(n,\lambda)}\right]}.$$
Let $g_{\tilde{D}-1}(x)$ and $g_{{\tilde{B}}_n-1}(x)$ be the probability generating functions of $\tilde{D}-1$ and $\tilde{B}_n-1$ respectively. Then, \eqref{eq.laplacebinomial} implies that
$$g_{{\tilde{B}}_n-1}(x)=g_{\tilde{D}-1}\left(1-p(\lambda,n)+p(\lambda,n)x\right).$$
Note that, for $g_{{\tilde{B}}_n-1}^{\circ k }$, the $k^{th}$ iterate of $g_{{\tilde{B}}_n-1}$, condition ({CM4}) is equivalent to
\begin{equation}\label{eq.treesdie}\liminf_{n\to\infty}\left(g_{{\tilde{B}}_n-1}^{\circ \lfloor \delta n^{(\alpha-1)/(\alpha+1)}\rfloor }(0)\right)^{\lfloor n^{1/(\alpha+1)}\rfloor}>0.\end{equation}
(See for instance the discussion below Theorem 2.3.1 in \cite{Duquesne2005}.)
As in the proof of  \cite[Proposition 5.25]{broutinLimitsMultiplicativeInhomogeneous2020}, it sufficient if we show that, for $t\in (0,\infty)$, and $r_n(t)$ the value such that
$$\int_{r_n(t)}^1 \frac{dr}{g_{{\tilde{B}}_n-1}(1-r)-1+r}=n^{(\alpha-1)/(\alpha+1)}t,$$
we have that
$$\limsup_{n\to\infty}n^{1/(\alpha+1)}r_n(t)<\infty.$$
Note that 
\begin{align*}
    \int_{r_n(t)}^1 \frac{dr}{g_{{\tilde{B}}_n-1}(1-r)-1+r}&=\int_{r_n(t)}^1 \frac{dr}{g_{\tilde{D}-1}\left(1-p(\lambda,n)r\right)-1+r}\\
    &=\int_{p(\lambda,n)r_n(t)}^{p(\lambda,n)} \frac{ds}{p(\lambda,n)g_{\tilde{D}-1}\left(1-s\right)-p(\lambda,n)+s},
\end{align*}
where we change variable to $s=p(\lambda,n)r$. Elementary calculation yields that 
$$g_{\tilde{D}-1}\left(1-s\right)=1-\left(\rho/\mu-1\right)s+\frac{C_\alpha}{\mu}s^\alpha+o(s^\alpha)$$ 
as $s\to 0$, which implies that
$$p(\lambda,n)g_{\tilde{D}-1}\left(1-s\right)-p(\lambda,n)+s=-\lambda n^{-(\alpha-1)/(\alpha+1)}s+p(\lambda,n)\frac{C_\alpha}{\mu}s^\alpha+p(\lambda,n)o(s^\alpha)$$
as $s\to 0$. Then, setting $v=n^{1/(\alpha+1)}s$ implies that
\begin{align*}
&\int_{r_n(t)}^1 \frac{dr}{g_{{\tilde{B}}_n-1}(1-r)-1+r}=n^{(\alpha-1)/(\alpha+1)}\int_{n^{1/(\alpha+1)}p(\lambda,n)r_n(t)}^{n^{1/(\alpha+1)}p(\lambda,n)} \frac{dv}{-\lambda v +\frac{C_\alpha}{\mu}v^\alpha+o(v^\alpha)},
\end{align*}
so that, since $p(n,\lambda)=\Theta(1)$ as $n\to \infty$, we obtain
$$\limsup_{n\to\infty}n^{1/(\alpha+1)}r_n(t)<\infty$$
as required. This implies that condition (CM4) is satisfied and so, indeed, performing bond percolation on the half-edges of a sequence of vertices with i.i.d.\ degrees in the supercritical regime yields a degree distribution that satisfies the conditions of Theorem \ref{thm.convgraphscriticalwindow}

\subsubsection{The methods in \cite{Conchon2018}}
In this section, we will further discuss the results and methods in  \cite{Conchon2018}. As mentioned previously, Conchon-Kerjan and Goldschmidt study a specific case of the model defined in Section \ref{subsec.model}, namely the case where the degree sequence does not depend on $n$ and $\lambda=0$. They prove Theorem \ref{thm.convgraphscriticalwindow} for that family of models, which is the content of \cite[Theorem 1.1]{Conchon2018}. (Their result includes the case $\alpha=2$, which we do not consider here.) 

The limit object is referred to as the \emph{$\alpha$-stable graph} for $\alpha\in (1,2)$ (and the \emph{Brownian graph} for $\alpha=2$). They obtain an additional result that identifies the components of the limit object as $\R$-trees encoded by tilted excursions of an $\alpha$-stable spectrally positive L\'evy process for $\alpha\in (1,2)$ (and tilted excursions of a Brownian motion for $\alpha=2$) with additional vertex identifications.  We cannot obtain such a description of the limit components in Theorem \ref{thm.convgraphscriticalwindow} because of their lack of self-similarity. 

We now give an informal overview of the proof of \cite[Theorem 1.1]{Conchon2018}. Much of the proof transfers over to our setting without change, so after this overview we will focus on the parts of the proof that are different in our setting. The method in \cite{Conchon2018} relies on using the configuration model in a depth-first manner, which we describe below. 

From the description of the configuration model, it is clear that we can pick an order of connecting half-edges to our convenience. Hence, we will choose an order that makes it similar to a depth-first exploration process. First, sample an i.i.d.\ degree sequence $D_1,\dots, D_n$ with $D_1\geq 1$ almost surely. Start from a vertex $v$ chosen with probability proportional to $D_v$ and label its half-edges in an arbitrary way. We maintain a \emph{stack}, which is an ordered list of the half-edges that we have seen but have not yet explored. Add all the half-edges of $v$ to the stack, ordered according to their labels, with the lowest label being on top of the stack. From now on, at each step, if the stack is non-empty, take the half-edge from the top of the stack, and sample its pair uniformly among the unpaired half-edges, i.e.\ the remaining half-edges on the stack, and the unexplored half-edges not on the stack. If the paired half-edge was not on the stack, say it was linked to vertex $w$, remove the paired half-edges from the system and place the remaining $D_w-1$ half-edges of $w$ on the top of the stack, arbitrarily labelled and in decreasing order of label, such that the lowest label of a half-edge of $w$ is now on top of the stack (unless the degree of $w$ is $1$). If the paired half-edge was on the stack, remove both paired half-edges from the system. If the stack is empty, we start a new connected component by selecting an unexplored vertex with probability proportional to its degree, and putting its half-edges on top of the stack. 

The argument then proceeds as follows.

\begin{enumerate}
    \item Conditionally on $D_1,\dots,D_n$, if we order the vertices by the time their first half-edge is paired in the configuration model, the ordered degree sequence $\hat{D}^n_1,\dots, \hat{D}^n_n$ is a size-biased random ordering of $D_1,\dots,D_n$, and the forest encoded by the \L ukasiewicz path $\tilde{S}^n_m=\sum_{i=1}^m\left(\hat{D}^n_i-2\right)$ is closely related to the components of the multigraph given by the configuration model. 
    \item For $i\neq j$, in general,  $\hat{D}^n_i$ is not equal to $\hat{D}^n_j$ in distribution, and  $\hat{D}^n_i$ and $\hat{D}^n_j$ are dependent. These facts makes $\tilde{S}^n$ hard to study. However, for $Z_1,\dots, Z_n$ i.i.d.\ with $\P(Z_1=k)=\frac{k\P(D_1=k)}{\E[D_1]},$ for $m\leq n$ there exists a function $\phi^n_m$ such that for $g:\N^m\to\R$, 
    $$\E[g(\hat{D}_1^n,\hat{D}_2^n,\dots,\hat{D}_m^n)]=\E[\phi_m^n(Z_1,Z_2,\dots,Z_m)g(Z_1,Z_2,\dots,Z_m)].$$ Moreover, $\phi^n_m$ behaves well under rescaling, which allows the authors of \cite{Conchon2018} to study $S_m:=\sum_{i=1}^m\left(Z_i-2\right)$, and then use the measure change to translate results on this simpler process to results on $\tilde{S}^n$. 
    \item Indeed, under rescaling, $S$ converges to an $\alpha$-stable spectrally positive L\'evy process $L$, jointly with its height process, and this result is used to show that $\tilde{S}$ (up to time $O(n^{\alpha/(\alpha+1)})$) converges to a process that is locally absolutely continuous to $L$, jointly with its height process.
    \item The excursions of $\tilde{S}^n$ above its running infimum and the corresponding excursions above $0$ of its height process encode individual trees in the forest. It is shown that the longest excursions explored up to time $\Theta(n^{\alpha/(\alpha+1)})$ converge under rescaling. The theory of size-biased point processes, developed by Aldous in \cite{Aldous1991}, is then used to show that, in fact, with high probability, all large excursions are observed in the first $\Theta(n^{\alpha/(\alpha+1)})$ steps, and that the excursions listed in decreasing order of length converge as well.
    \item By adding extra randomness and making some vertex identifications, the forest encoded by $\tilde{S}^n$ can be modified to yield a multigraph that is equal in law to the graph created by the configuration model, and these modifications behave well under rescaling the graph and taking limits.
    \item Finally, the authors show that conditioning on the graph not containing multiple edges and loops does not affect the distribution of the largest components. This follows by adapting an argument of Joseph in \cite{Joseph2014}, which shows that the first loops and multiple edges are sampled far beyond the time scale $\Theta(n^{\alpha/(\alpha+1)})$, and so their presence or absence cannot affect the scaling limit. 
\end{enumerate}

\subsection{Adapting the methods in \cite{Conchon2018} to the critical window}
The largest barrier to generalising the methods in \cite{Conchon2018} is showing the convergence under rescaling of $S$, jointly with its height process. The results proved in Section \ref{sec.convergence} allow us to do this. After that, it is trivial to extend most of the arguments in \cite{Conchon2018} to the critical window. 

The convergence under rescaling of $S$ is the content of Proposition \ref{prop.scalingcriticalwindow}. Then, we discuss the results in \cite{Conchon2018} that are not
trivially extended to the critical window and need some further justification. 

Let $D_1^n,\dots, D_n^n$ be i.i.d.\ with a degree distribution as specified in \ref{subsec.model}. Recall that the degree distribution depends on both $\lambda$ and $n$, but that we have dropped the dependency on $\lambda$ in the notation. 

We consider the configuration model executed in depth first order on vertices with degrees $D_1^n,\dots, D_n^n$.  Let $(\hat{D}_1^n,\dots,\hat{D}_n^n)$ denote the degrees in order of discovery, such that $(\hat{D}_1^n,\dots,\hat{D}_n^n)$ is distributed as a size-biased random ordering of $D_1^n,\dots,D_n^n$. This is defined as follows.

Let $\Sigma$ be a random permutation of $\{1,\dots,n\}$ such that
$$\P(\Sigma=\sigma|D_1^n,\dots,D_n^n)=\frac{D^n_{\sigma(1)}}{\sum_{j=1}^n D^n_{\sigma(j)}}\frac{D_{\sigma(2)}}{\sum_{j=2}^n D^n_{\sigma(j)}}\cdots \frac{D^n_{\sigma(n)}}{D^n_{\sigma(n)}}.$$
Then, by Proposition 3.2 of \cite{Conchon2018}, 
\begin{align*}&\P (\hat{D}_1^n=k_1, \hat{D}_2^n=k_2,\dots, \hat{D}_n^n=k_n)\\
&=n! k_1 \P(D_1^n=k_1) k_2 \P(D_1^n=k_2) \cdots k_n \P(D_1^n=k_n) \prod\limits_{i=1}^n \frac{1}{\sum_{j=i}^n k_j}.\end{align*}
Now, let $0\leq m\leq n$, and $k_1,k_2,\dots,k_m\geq 1$, and define $\Xi^n_{n-m}$ to be a random variable having the same law as $D^n_{m+1}+\cdots+D^n_{n}$. Then, for 
$$\phi_m^n(k_1,k_2,\dots,k_m):=\frac{n!\mu^m}{(n-m)!}\E\left[\prod \limits_{i=1}^m \frac{1}{\sum_{j=i}^m k_j+\Xi^n_{n-m}}\right],$$
Proposition 3.2 in \cite{Conchon2018} yields that for $Z^n_1,\dots,Z^n_n$ i.i.d.\ random variables with the size-biased degree of $D_1^n$, i.e.
$$\P(Z^n_1=k)=\frac{k\P(D^n_1=k)}{\E[D^n_1]},$$
for any test-function $g:\N^m\rightarrow \R,$
$$\E[g(\hat{D}_1^n,\hat{D}_2^n,\dots,\hat{D}_m^n)]=\E[\phi_m^n(Z^n_1,Z^n_2,\dots,Z^n_m)g(Z^n_1,Z^n_2,\dots,Z^n_m)],$$
i.e., $\phi_m^n$ defines a measure change to get from a vector of size-biased distributed random variables to a vector of size-biased ordered random variables. We note that $$\tilde{S}^n(k):=\sum_{i=1}^k\left(\hat{D}^n_i-2\right)$$
is the \L ukasiewicz path of a forest that is closely related to the depth-first spanning forest of our graph of interest, because it encodes the degrees in order of discovery.  Therefore, the existence of the measure change motivates the study of the limit under rescaling of $${S}^n(k):=\sum_{i=1}^k\left(Z^n_i-2\right)$$
and its corresponding height process. This is the content of the following proposition.

\begin{proposition}\label{prop.scalingcriticalwindow}
Let $L$ be a spectrally positive $\alpha$-stable L\'evy process with L\'evy measure $\pi(dx)=\frac{c}{\mu}x^{-(\alpha+1)}dx$. Then, for any $\lambda\in \R$ there exists a continuous modification of the height process of $(L_t+\lambda t, t\geq 0)$, which we will denote by $H^\lambda$. Moreover, for $H^n$ the height process corresponding to $S^n$, we have that, as $n\rightarrow \infty$,
\begin{equation*}\left(  n^{-1/(\alpha+1)}S^n\left(\lfloor tn^{{\alpha/(\alpha+1)}}\rfloor\right), n^{-(\alpha-1)/(\alpha+1)}H^n\left(\lfloor tn^{{\alpha/(\alpha+1)}}\rfloor\right) t\geq 0\right)\xrightarrow{d} (L_t+\lambda t, H^\lambda_t, t\geq 0)\end{equation*}
in $\D(\R_+,\R)^2$.
\end{proposition}

We will use Theorem \ref{thm.mainconvergenceresult}  to prove Proposition \ref{prop.scalingcriticalwindow}. We will first study the Laplace transform of $Z_1^n$, which is the content of the following lemma. 

\begin{lemma}\label{lemma.laplaceoperatorZn}
Define $\cL_{Z_1^{n}}(s)=\E[\exp(-s Z_1^{n})] $. Then, 
$$\cL_{Z_1^{n}}(s)=1-\left(2+\lambda n ^{-{(\alpha-1)/(\alpha+1)}}\right)s+\frac{C^{(n)}}{\mu_n}s^{\alpha}+o(s^{\alpha})$$
as $s\to 0$.
\end{lemma}
\begin{proof}
Note that 
\begin{align}\begin{split}\label{eq.Laplace transformnoscaling}
\cL_{Z_1^{n}}''(s)&=\E \left[(Z_1^{n})^2\exp(-s Z_1^{n})\right]\\
&=\sum\limits_{k=1}^{\infty}\frac{k^3\nu_k^{n}e^{-sk}}{\mu_n}\\
&=\frac{c_n}{\mu_n}s^{\alpha-2}\Gamma(2-\alpha)+o(s^{\alpha-2})
\end{split}\end{align}
for $s\rightarrow 0$, where the last equality follows from the Euler-Maclaurin formula, using that $\P( Z_1^{n}= k)\sim \frac{c_n}{\mu_n}k^{-(\alpha+1)}$ as $k\rightarrow \infty$.
Then, because $\E[Z_1^{n} ]=2+\lambda n ^{-{(\alpha-1)/(\alpha+1)}}$, integrating twice gives the result.
\end{proof}


\begin{proof}[Proof of Proposition \ref{prop.scalingcriticalwindow}]
We will first prove that $S^n$ converges in distribution under rescaling. Set $$M^{n}(k)=\sum\limits_{i=1}^{k}\left(Z_i^{n}-2-\lambda n ^{-{(\alpha-1)/(\alpha+1)}}\right)$$ and $$A^{n}(k)= k\lambda n ^{-{(\alpha-1)/(\alpha+1)}}$$ so that $$S^n(k)=M^{n}(k)+A^{n}(k)$$ 
is the Doob-Meyer decomposition of $S^n$.
Firstly, observe that $$\left(n^{-1/(\alpha+1)} A^{n}\left(\lfloor tn^{{\alpha/(\alpha+1)}}\rfloor\right),t \geq 0\right) \rightarrow (\lambda t,t \geq 0) $$
in $\D(\R_+,\R)$  as $n\to\infty$.

We claim that for every $t\geq 0$,
\begin{equation}\label{eq.claimMn} n^{-1/(\alpha+1)} M^{n}\left(\lfloor tn^{{\alpha/(\alpha+1)}}\rfloor\right)\xrightarrow{d} L_t .\end{equation} Firstly, observe that for all $\theta>0$,
$$\E\left[\exp(-\theta L_t)\right]=\exp\left(t\frac{c\Gamma(2-\alpha)}{\mu \alpha (\alpha-1)}\theta^\alpha\right).$$
Recall that $C_{\alpha}=\frac{c\Gamma(2-\alpha)}{\alpha(\alpha-1)}$, and set $C^{(n)}=\frac{c_n\Gamma(2-\alpha)}{\alpha(\alpha-1)}$, so that $C^{(n)}\to C_\alpha$ as $n\to\infty$.
We will show that for every $\theta>0$,
$$\E\left[\exp\left(-\theta n^{-1/(\alpha+1)} M^{n}\left(\lfloor tn^{{\alpha/(\alpha+1)}}\rfloor\right) \right)\right] \rightarrow \exp\left(t\frac{C_\alpha}{\mu }\theta^\alpha\right)$$ 
as $n\rightarrow \infty$, which will prove the claim.

Note that \begin{align}\begin{split}\label{eq.laplace} &\E\left[\exp(-\theta n^{-1/(\alpha+1)} M^{n}(\lfloor tn^{{\alpha/(\alpha+1)}}\rfloor)\right ]\\&=\left[\cL_{Z_1^{n}}(\theta n^{-1/(\alpha+1)})\exp\left(\theta n^{-1/(\alpha+1)}\left(2+\lambda n ^{-{(\alpha-1)/(\alpha+1)}}\right)\right) \right]^{\lfloor tn^{{\alpha/(\alpha+1)}}\rfloor }.\end{split}\end{align}
Then, Lemma \ref{lemma.laplaceoperatorZn} implies that
\begin{align*} \cL_{Z_1^{n}}(s)=\exp\left(-\left(2+\lambda n ^{-{(\alpha-1)/(\alpha+1)}}\right)s+\frac{C^{(n)}}{\mu_n}s^{\alpha}+o(s^{\alpha})\right)
\end{align*}
as $s\rightarrow 0$. Plugging this into \eqref{eq.laplace}, we find that as $n\rightarrow \infty$,
$$\E\left[\exp\left(-\theta n^{-1/(\alpha+1)} M^{n}(\lfloor tn^{{\alpha/(\alpha+1)}}\rfloor)\right)\right ]\rightarrow \exp\left(t\frac{C_{\alpha}}{\mu }\theta^\alpha\right),$$
which proves \eqref{eq.claimMn}.

Now, using \cite[Theorem 16.14]{Kallenberg2002}, we may deduce that as $n\rightarrow \infty$,
\begin{equation}\label{eq.convergenceS} \left( n^{-1/(\alpha+1)} S^n\left(\lfloor tn^{{\alpha/(\alpha+1)}}\rfloor\right), t\geq 0\right)\xrightarrow{d} (L_t+\lambda t, t\geq 0)\end{equation}
in $\D(\R_+,\R)$.

In order to also obtain the convergence of the height process under rescaling, we note that for $\lambda\leq 0$, we can directly apply \cite[Theorem 1.4.3 and Theorem 2.2.1]{Duquesne2005}, stated in this work as Theorem \ref{thm.convergencesubcritical}. In the case of $\lambda>0$, we apply Theorem \ref{thm.mainconvergenceresult}. In both results, we use the scaling parameters $p=n^{1/\alpha+1}$ and $\gamma_p=n^{(\alpha-1)/(\alpha+1)}$. The conditions of the theorems then follow directly from our assumptions on the degree distributions and \eqref{eq.convergenceS}. This implies the existence of a continuous modification of the height process and the claimed convergence. 
\end{proof}

We will now prove that in the cases we consider, $\phi^n_m$ behaves well under rescaling. This is the content of the following proposition, which is a generalization of the proof of \cite[Proposition 4.3]{Conchon2018}. 

\begin{proposition}\label{prop.convergencemeasurechange}
Set 
$$\phi(t)=\exp\left(\frac{1}{\mu} \int_0^t s dL_s -\frac{C_{\alpha}}{(\alpha+1) \mu^{\alpha+1}}t^{\alpha+1} \right). $$
Then 
$$\phi(n,t)\xrightarrow{d}\phi(t)$$ as $n\rightarrow \infty$, and the sequence $\{\phi(n,t), n\in\N\}$ is uniformly integrable.
\end{proposition}
The proof will follow the structure of the proof of \cite[Proposition 4.3]{Conchon2018}, but we will need to adapt the technical lemmas presented there to our more general setting.

\begin{proof}
The following technical lemmas need justification in the more general setting. 
\begin{itemize}
    \item For $\lambda=0$ and a degree distribution that does not depend on $n$, for any $t\geq 0$, by \cite[Lemma 4.4]{Conchon2018} it holds that, as $n\to\infty$, $$\frac{1}{n}\sum\limits_{k=0}^{\lfloor tn^{{\alpha/(\alpha+1)}}\rfloor -1} S^n(k) \xrightarrow{d} \int_0^tL_sds.$$
    Using the continuous mapping theorem (see e.g.\ \cite[Theorem 3.2.4]{Durrett2010}), we find that for general $\lambda$ and for the degree distribution depending on $n$, for any $t\geq 0$,
\begin{equation}\label{eq.integral}\frac{1}{n}\sum\limits_{k=0}^{\lfloor tn^{{\alpha/(\alpha+1)}}\rfloor -1} S^n(k) \xrightarrow{d} \int_0^tL_sds+\frac{\lambda t^2}{2}.\end{equation}
\item Define $\cL_{D_1^n}(s):=\E[\exp(-sD_1^{n})]$. Then, \cite[Lemma 4.5a]{Conchon2018} and the argument thereafter state that for $\lambda=0$ and a degree distribution that does not depend on $n$,
\begin{align*}
\cL_{D_1^n}(s)&=1-\mu_ns+\mu_ns^2-\frac{C_{\alpha}}{\alpha+1}s^{\alpha+1}+o(s^{\alpha+1})\\
&=\exp\left(-\mu_n s + \frac{\mu_n(2-\mu_n)}{2}s^2-\frac{C_{\alpha}}{\alpha+1}s^{\alpha+1}+o(s^{\alpha+1})\right)
\end{align*}
In our set-up, we find, similarly to how we obtained $\cK_{n}''(s)$, that 
$$\cL_{D_1^n}'''(s)=-c\Gamma(2-\alpha)s^{\alpha-2}+o(s^{\alpha-2})$$
as $s\rightarrow 0$, and using $\E[D_1^{n}]=\mu_n$ and $\E \left[(D_1^{n})^2\right]=\left(2+\lambda n ^{-{(\alpha-1)/(\alpha+1)}}\right) \mu_n$, we get, by integrating three times, that
\begin{align}
\begin{split}
\cL_{D_1^n}(s)&=1-\mu_ns+\left(1+\frac{\lambda}{2} n^{-{(\alpha-1)/(\alpha+1)}}\right)\mu_ns^2-\frac{C^{(n)}}{\alpha+1}s^{\alpha+1}+o(s^{\alpha+1})\\
&=\exp\left(-\mu_n s + \frac{\mu_n(2+\lambda n ^{-{(\alpha-1)/(\alpha+1)}}-\mu_n)}{2}s^2-\frac{C^{(n)}}{\alpha+1}s^{\alpha+1}+o(s^{\alpha+1})\right)\label{eq.Laplace transform}
\end{split}
\end{align}
as $s\rightarrow 0$.
\item  For $\lambda=0$ and a degree distribution that does not depend on $n$, for $m=\lfloor tn^{\alpha/(\alpha+1)}\rfloor$, by \cite[Lemma 4.6]{Conchon2018} it holds that, as $n\to\infty$, 
\begin{equation*}\exp \left(m- \frac{ 2+\mu_n}{2\mu_n}\frac{m^2}{n}\right)\left[\cL\left(\frac{m}{n\mu_n}\right)\right] ^{n-m} \rightarrow \exp\left(- \frac{C_{\alpha}}{(\alpha+1 )\mu^{\alpha+1}}t^{\alpha+1}\right).\end{equation*}
By \eqref{eq.Laplace transform} we straightforwardly obtain that, in our set-up, as $n\rightarrow \infty$, 
\begin{equation}\label{eq.exponlimit}\exp \left(m- \frac{ 2+\mu_n}{2\mu_n}\frac{m^2}{n}\right)\left[\cL\left(\frac{m}{n\mu_n}\right)\right] ^{n-m} \rightarrow \exp\left(- \frac{C_{\alpha}}{(\alpha+1 )\mu^{\alpha+1}}t^{\alpha+1}+\frac{\lambda t^2}{2\mu}\right).\end{equation}
\item For $\lambda=0$ and a degree distribution that does not depend on $n$, for $s(0)=0$ and $s(i)=\sum_{j=1}^i (k_j-2)$, for $i\geq 1$, by \cite[Lemma 4.7]{Conchon2018},  it holds that $$\phi_m^{n}(k_1,\dots,k_m)\geq \exp\left(\frac{1}{n\mu_n}\sum\limits_{i=0}^{m}(s(i)-s(m))-\frac{C_\alpha}{(\alpha+1 )\mu^{\alpha+1}}t^{\alpha+1} \right)(1+o(1)) .$$
Adapting the proof to our set-up, using  \eqref{eq.exponlimit}, we find that for $s(0)=0$ and $s(i)=\sum_{j=1}^i (k_j-2)$, for $i\geq 1$, 
$$\phi_m^{n}(k_1,\dots,k_m)\geq \exp\left(\frac{1}{n\mu_n}\sum\limits_{i=0}^{m}(s(i)-s(m))-\frac{C_\alpha}{(\alpha+1 )\mu^{\alpha+1}}t^{\alpha+1}+\frac{\lambda t^2}{2\mu}\right)(1+o(1)) .$$ The proof of this can be found in the Appendix. 
\end{itemize}

Now, we have that $$\phi(n,t)\geq  \underline{\phi}(n,t):= \exp\left(\frac{1}{n\mu_n}\sum\limits_{i=0}^{m}(S^n(i)-S^n(m))-\frac{C_\alpha}{(\alpha+1) \mu^{\alpha+1}}t^{\alpha+1}+\frac{\lambda t^2}{2\mu}\right).$$

But then, by \eqref{eq.convergenceS} and \eqref{eq.integral}, we get that, 
since 
$$\frac{1}{n}\sum\limits_{k=0}^{\lfloor tn^{{\alpha/(\alpha+1)}}\rfloor -1} \left(S^n(k)-S^n(\lfloor tn^{{\alpha/(\alpha+1)}}\rfloor -1)\right) \xrightarrow{d} \int_0^t(L_s-L_t) ds-\frac{\lambda t^2}{2},$$

$$ \underline{\phi}(n,t)\xrightarrow{d}\exp \left(\frac{1}{\mu}\int_0^t (L_s-L_t)ds-\frac{C_{\alpha}}{(\alpha+1) \mu^{\alpha+1}}t^{\alpha+1}\right).$$
Then, we finish the proof like the proof of \cite[Proposition 3.3]{Conchon2018} to obtain
$$\phi(n,t)\xrightarrow{d}\phi(t)$$ as $n\rightarrow \infty$, and that the sequence $\{\phi(n,t), n\in\N\}$ is uniformly integrable. The details can be found in the Appendix.
\end{proof}

Remember that  $n^{-1/(\alpha+1)}\left( S^n\left(\lfloor tn^{{\alpha/(\alpha+1)}}\rfloor\right), t\geq 0\right)$ converges in law to $(L_t+\lambda t, t\geq 0)$ as $n\to \infty$. Via the measure change $\phi(n,t)$ we can get from $\left(S^n\left(\lfloor sn^{{\alpha/(\alpha+1)}}\rfloor\right),0\leq s \leq t\right)$ to $\left(\tilde{S}^{n}\left(\lfloor sn^{{\alpha/(\alpha+1)}}\rfloor\right),0\leq s \leq t \right)$. The random variable $\phi(n,t)$ converges in law to $\phi(t)$ as $n\to\infty$. Therefore, we will define the process $(\Tilde{K}^{\lambda},\Tilde{H}^\lambda)$ via the following measure change. For $t\geq 0$, for every non-negative integrable functional $F:\D([0,t],\R^2)\to\R$, 
\begin{align}\begin{split}\label{def.tildeK} \E[F(\Tilde{K}^{\lambda}_s,\Tilde{H}^{\lambda}_s, 0\leq s \leq t)]&=\E\left[\phi(t)F(L_s+\lambda s,H^\lambda_s,0\leq s \leq t)\right].
\end{split}\end{align}

\begin{proposition}\label{prop.jointconvergencesizebiasedordered}
We have
\begin{equation*} \left(n^{-1/(\alpha+1)} \tilde{S}^{n}\left(\lfloor sn^{{\alpha/(\alpha+1)}}\rfloor\right), n^{-\alpha/(\alpha+1)} \tilde{H}^{n}\left(\lfloor sn^{{\alpha/(\alpha+1)}}\rfloor\right), 0\leq s \leq t\right)\xrightarrow{d} (\tilde{K}^\lambda_s, \tilde{H}^\lambda_s ,0\leq s \leq t)\end{equation*}
in $\D(\R_+,\R^2)$ as $n\to\infty$ . 
\end{proposition}
\begin{proof}

We want to show that for any $t\geq 0$ and any bounded continuous test function $f:\D([0,t],\R^2) \rightarrow \R$, for $\tilde{H}^n$ the height process corresponding to $\tilde{S}^n$,
\begin{align*}&\E\left[f\left( n^{-1/(\alpha+1)} \tilde{S}^{n}\left(\lfloor sn^{{\alpha/(\alpha+1)}}\rfloor\right), n^{-\alpha/(\alpha+1)} \tilde{H}^{n}\left(\lfloor sn^{{\alpha/(\alpha+1)}}\rfloor\right), 0\leq s \leq t\right)\right] \\
&\rightarrow \E \left[f(\tilde{K}^\lambda_s, \tilde{H}^\lambda_s ,0\leq s \leq t)\right],\end{align*}
as $n\rightarrow \infty$. 
By using our measure change, this is equivalent to showing that for any $t\geq 0$ and any bounded continuous test function $f:\D([0,t],\R^2)\rightarrow \R$,
\begin{align*}&\E \left[\phi(n,t) f\left( n^{-1/(\alpha+1)} {S}^{n}\left(\lfloor sn^{{\alpha/(\alpha+1)}}\rfloor\right), n^{-\alpha/(\alpha+1)}{H}^{n}\left(\lfloor sn^{{\alpha/(\alpha+1)}}\rfloor\right), 0\leq s \leq t\right) \right] \\
&\rightarrow \E [\phi(t)f(L_s+\lambda s,H^\lambda_s, 0\leq s\leq t)].\end{align*}
We now finish as in the proof of \cite[Theorem 4.1]{Conchon2018} in order to obtain the desired result. \end{proof}

The following proposition characterises the law of $\tilde{K}^\lambda$.
\begin{proposition}\label{prop.laplaciantildeK}
For $\tilde{K}^\lambda(t)$ as defined in \eqref{def.tildeK}, for $\theta>0$
$$\E\left[\exp\left(-\theta \tilde{K}^\lambda(t) \right)\right]=\exp\left( -\theta \lambda t + \theta C_\alpha \frac{t^\alpha}{\mu^\alpha}+\int_0^t ds \int_0 ^\infty \frac{c}{\mu}x^{-(\alpha+1)}e^{-xs/\mu}dx(e^{-\theta x}-1+\theta x)\right) .$$
\end{proposition}
The proof can be found in the appendix.

\begin{proof}[Proof of Theorem \ref{thm.convgraphscriticalwindow}] Given Proposition \ref{prop.jointconvergencesizebiasedordered}, the rest of the proof is completely analogous to the proof of in \cite[Theorem 1.1]{Conchon2018}. Proposition \ref{prop.laplaciantildeK} characterizes the law of the encoding process of the components of the limit object. 
\end{proof}

\section*{Acknowledgements}
The author would like to thank her advisor Christina Goldschmidt for many productive meetings. Moreover, she would like to thank Matthias Winkel for helpful input in the early stages of this project. She would like to thank Thomas Duquesne for pointing out interesting literature. Finally, she would like to thank Thomas Hughes for careful proofreading. 
\appendix
\section{Proofs of technical results}\label{app.proofs}
\begin{proof}[Proof of Lemma \ref{lemma.technicalconvergenceinfimum}]
Firstly, note that since $f$ and $g$ are increasing, also the function given by $\left(\inf \{s\geq 0:f(t)<g(s)\},t \geq 0 \right)$ is increasing, and so in particular it has limits from the left and from the right at every point of its domain. Fix $t\geq 0$. Suppose that
$$\lim\limits_{q\uparrow t} \inf\{s\geq 0:f(q)<g(s)\}<\inf\{s\geq 0:f(t)<g(s)\}.$$
Then we must have that
\begin{enumerate}
    \item there is an $\tilde{s}$ such that $f(t)=g(\tilde{s})$, and
    \item $f(q)<f(t)$ for all $q<t$.
\end{enumerate}
It follows that $$\lim\limits_{q\uparrow t} \inf\{s\geq 0:f(q)<g(s)\}=\tilde{s}.$$
But then, since $g$ is strictly increasing, 
$$\inf \{s\geq 0:f(t)<g(s)\}=\inf \{s\geq 0:g(\tilde{s})<g(s)\}=\tilde{s},$$
so we must have
$$\lim\limits_{q\uparrow t} \inf\{s\geq 0:f(q)<g(s)\}=\inf\{s\geq 0:f(t)<g(s)\}.$$
We also need to show that
$$\lim\limits_{q\downarrow t} \inf\{s\geq 0:f(q)<g(s)\}=\inf\{s\geq 0:f(t)<g(s)\}.$$
Fix $\epsilon>0$. Suppose $\inf\{s\geq 0:f(t)<g(s)\}=\tilde{s}$, which implies that $f(t)<g(\tilde{s}+\epsilon)$. Then, as $q\downarrow t$, we have $f(q)\downarrow f(t)$, so for $q>t$ small enough, $f(q)<g(\tilde{s}+\epsilon)$. Hence,
$$\inf\{s\geq 0:f(q)<g(s)\}\leq \tilde{s}+\epsilon ,$$
which proves that 
$$\lim\limits_{q\downarrow t} \inf\{s\geq 0:f(q)<g(s)\}=\inf\{s\geq 0:f(t)<g(s)\}.$$
Therefore, $\left(\inf \{s\geq 0:f(t)<g(s)\},t \geq 0 \right)$ is continuous.\\
By Ethier and Kurtz \cite[Proposition 3.6.5]{Ethier1986}, proving
$$\left(\inf \{s\geq 0:f_n(t)<g_n(s)\},t \geq 0 \right)\to \left(\inf \{s\geq 0:f(t)<g(s)\},t \geq 0 \right)$$
in $\D(\R_+,\R)$  as $n\to\infty$ is then equivalent to showing that for all $t\in \R_+$, and for all $(t_n)_{n\geq 0}$ in $\R_+$ such that $t_n\to t$, we have 
$$\inf\{s\geq 0:f_n(t_n)<g_n(s)\}\to \inf\{s\geq 0:f(t)<g(s)\}.$$
 Suppose $\inf\{s\geq 0:f(t)<g(s)\}=\tilde{s}$. Fix $\epsilon>0$ and $(t_n)_{n\geq 0}$ in $\R_+$ such that $t_n\to t$. \\
We will first show that for $n$ large enough, $$\inf\{s\geq 0:f_n(t_n)<g_n(s)\}\leq \tilde{s}+\epsilon.$$
By definition of $\tilde{s}$ and monotonicity, we have that $f(t)<g(\tilde{s}+\epsilon/2)$. Since $g$ is strictly increasing, we have that there is a $\delta>0$ such that $g(\tilde{s}+\epsilon/2)=f(t)+\delta$. Moreover, since $f_n\to f$ in $\D(\R_+,\R)$ as $n\to\infty$ and $f$ is continuous, $f_n(t_n)\to f(t)$, so we may pick $n$ large enough such that 
$$|f_n(t_n)-f(t)|<\delta/2.$$
Since $g_n\to g$, for $n$ large enough there exists a monotone bijection $\lambda_n:[0,\tilde{s}+2\epsilon]\to [0,\tilde{s}+2\epsilon]$ such that
$$\sup_{t\in [0,\tilde{s}+2\epsilon]}|\lambda_n(t)-t|<\epsilon/4$$
and 
$$\sup_{t\in [0,\tilde{s}+2\epsilon]} |g_n(\lambda_n(t))-g(t)|<\delta/2.$$
Hence, 
$$f_n(t_n)<f(t)+\delta/2< g_n(\lambda_n(\tilde{s}+\epsilon/2))$$
and 
$$f_n(t_n)<g_n(s)$$
for $s=\tilde{s}+3\epsilon/4$. Hence, for $n$ large enough,
$$\inf\{s\geq 0:f_n(t_n)<g_n(s)\}<\tilde{s}+\epsilon.$$
We now want to show that for $n$ large enough,
$$\inf\{s\geq 0:f_n(t_n)<g_n(s)\}\geq \tilde{s} -\epsilon.$$
Fix $s<\tilde{s}-\epsilon$. We will be done if we can show that $f_n(t_n)\geq g_n(s)$. 
Since $g$ is strictly increasing, we know that there exists a $\delta>0$ such that $g(\tilde{s}-\epsilon/2)+\delta =g(\tilde{s}-\epsilon/4)$. Also, by definition of $\tilde{s}$, we have $f(t)\geq g(\tilde{s}-\epsilon/4)$.
Pick $n$ large enough that there exists a monotone bijection $\lambda_n:[0,\tilde{s}+2\epsilon]\to [0,\tilde{s}+2\epsilon]$ such that
$$\sup_{t\in [0,\tilde{s}+2\epsilon]} |\lambda_n(t)-t|<\epsilon/4$$
and 
$$\sup_{t\in [0,\tilde{s}+2\epsilon]} |g_n(\lambda_n(t))-g(t)|<\delta/2,$$
and such that
$$|f_n(t_n)-f(t)|<\delta/2.$$
Then, 
\begin{align*} g_n(s)&\leq g_n(\tilde{s}-\epsilon)\\
&\leq g(\tilde{s}-\epsilon/2)+\delta/2\\
&=g(\tilde{s}-\epsilon/4)-\delta/2\\
&\leq f(t)-\delta/2\\
&\leq f_n(t_n),
\end{align*}
which proves the statement.\\
Finally, to show that for $(\epsilon_n)_{n\geq 0}$ such that $\epsilon_n\downarrow 0$,
$$\left(\inf \{s\geq 0:f_n(t-\epsilon_n s)<g_n(s)\},t \geq 0 \right)\to \left(\inf \{s\geq 0:f(t)<g(s)\},t \geq 0 \right)$$
in $\D(\R_+,\R)$  as $n\to\infty$,
fix $t\in \R_+$, $t_n\to t$ and $\epsilon>0$, and we need to show that
$$\inf \{s\geq 0:f_n(t_n-\epsilon_n s)<g_n(s)\}\to \inf \{s\geq 0:f(t)<g(s)\}=:\tilde{s}.$$
Firstly, note that 
$$\inf \{s\geq 0:f_n(t_n-\epsilon_n s)<g_n(s)\}\leq \inf \{s\geq 0:f_n(t_n)<g_n(s)\},$$
which is smaller that $\tilde{s}+\epsilon$ for $n$ large enough by previous results. Moreover, note that for $s\leq \tilde{s}-\epsilon$,
$$f_n(t_n-\epsilon_n s)\geq f_n(t_n-\epsilon_n (\tilde{s}-\epsilon)),$$
which is larger than $g_n(s)$ for $n$ large enough by previous results, since $t_n-\epsilon_n(\tilde{s}-\epsilon)\to t $. Hence, $$\inf \{s\geq 0:f_n(t_n-\epsilon_n s)<g_n(s)\}\geq \tilde{s}-\epsilon$$ for $n$ large enough, which concludes the proof. \end{proof}

We will now prove two technical lemmas needed for the proof of Proposition \ref{prop.convergencemeasurechange}.

\begin{lemma}For $m=\lfloor t n^{\alpha/(\alpha+1)}\rfloor$, for
$$\phi_m^{n}(k_1,\dots,k_m)= \frac{n!\mu^m_{n}}{(n-m)!}\E \left[\prod\limits_{i=1}^m \frac{1}{\sum_{j=i}^m k_j + \Xi_{n-m}}\right],$$
$s(0)=0$ and $s(i)=\sum_{j=1}^i (k_j-2)$, for $i\geq 1$, we have that
$$\phi_m^{n}(k_1,\dots,k_m)\geq \exp\left(\frac{1}{n\mu_{n}}\sum\limits_{i=0}^{m}(s(i)-s(m))-\frac{C^{(n)}}{(\alpha+1 )\mu^{\alpha+1}}t^{\alpha+1}+\frac{\lambda t^2}{2\mu}\right)(1+o(1)) .$$ \end{lemma}
The arguments presented are an adaption of the proof of \cite[Lemma 4.7]{Conchon2018}.

\begin{proof}
We can write 
\begin{align*}\phi_m^{n}(k_1,\dots,k_m)&=\prod\limits_{i=1}^{m-1}\left(1-\frac{i}{n}\right) \E \left[\prod\limits_{i=1}^{m}\left(\frac{n\mu_n}{\sum_{j=i}^m k_j+\Xi_{n-m}}\right)\right]\\
&=\E\left[\exp\left(\sum\limits_{i=1}^{m-1}\log\left(1-\frac{i}{n}\right)-\sum\limits_{i=1}^m \log \left(\frac{\Xi_{n-m}}{n\mu_n}+\frac{1}{n\mu_n}\sum\limits_{j=i}^m k_j\right)\right)\right].\end{align*}
Then, since for all $x\in (-1,\infty)$, $\log(1+x)\leq x$ and for all $1\leq i \leq m-1$, $\log(1-i/n)\geq -i/n+m^2/n^2$, we have that 
\begin{align*}
    \phi_m^{n}(k_1,\dots,k_m)\geq& \exp\left(\frac{-m(m-1)}{2n}-\frac{m^3}{n^2}-m+\frac{1}{n\mu_n}\sum\limits_{i=0}^m(s(i)-s(m))-\frac{m(m+1)}{n\mu_n}\right)\\
    &\times\E\left[\exp\left(\frac{-m}{n\mu_n}\Xi_{n-m}\right)\right]\\
    =&\exp\left(\frac{1}{n\mu_n}\sum\limits_{i=0}^m(s(i)-s(m))\right)\exp\left(m-\frac{(2+\mu_n)m^2}{2\mu_n n}\right)\\
    &\times \left[\cL_{D_1^n}\left(\frac{m}{n\mu_n}\right)\right]^{n-m} \exp\left(\frac{\mu_n-2}{2\mu_n n}-\frac{m^3}{n^2}\right)
\end{align*}
Note that by definition of $m$, $m^3/n^2=\Theta\left(n^{(\alpha-2)/(\alpha+1)}\right)=o(1)$, so the final exponent tends to $1$ as $n\to \infty$. 
Then by \eqref{eq.exponlimit}, which states that 
as $n\rightarrow \infty$, 
$$\exp \left(m- \frac{ 2+\mu_{n}}{2\mu_{n}}\frac{m^2}{n}\right)\left[\cL_{D_1^n}\left(\frac{m}{n\mu_{n}}\right)\right] ^{n-m} \rightarrow \exp\left(- \frac{C_{\alpha}}{(\alpha+1 )\mu^{\alpha+1}}t^{\alpha+1}+\frac{\lambda t^2}{2\mu}\right),$$ the desired result follows. 
\end{proof}
Furthermore, recall that $$\phi(n,t)=\phi_m^n(Z_1^n,\dots,Z_m^n)$$ and $$\phi(t)=\exp\left(-\frac{1}{\mu}\int_0^t s dL_s-\frac{C_\alpha}{(\alpha+1)\mu^{\alpha+1}} t^{\alpha+1}\right).$$
\begin{lemma} As $n\to \infty$,
$$\phi(n,t)\overset{d}{\to} \phi(t),$$
and  $\{\phi(n,t),n\in\N\}$ is uniformly integrable.
\end{lemma}
The arguments presented are an adaptation of the proof of \cite[Proposition 4.3]{Conchon2018}.
\begin{proof}
Recall that $S^n(k)=\sum_{i=1}^k(Z^n_i-2).$
By the statement proved above, we have that 
$$\phi(n,t)\geq \underline{\phi}(n,t):=\exp\left(\frac{1}{n\mu_{n}}\sum\limits_{i=0}^{m}(S^n(i)-S^n(m))-\frac{C^{(n)}}{(\alpha+1 )\mu^{\alpha+1}}t^{\alpha+1}+\frac{\lambda t^2}{2\mu}\right)(1+o(1)).$$
 Then, by Theorem \ref{prop.scalingcriticalwindow}
 and by \eqref{eq.integral},
 we get that 
$$\frac{1}{n\mu_{n}}\sum\limits_{i=0}^{m}(S^n(i)-S^n(m))\overset{d}{\to}\frac{1}{\mu}\int_0^t (L_s^\lambda -L_t^\lambda)ds-\frac{\lambda t^2}{2\mu}.$$
Hence, we get that 
\begin{align*}\underline{\phi}(n,t)\overset{d}{\to} &\exp\left(\frac{1}{\mu} \int_0^t (L_s-L_t)ds-\frac{C_\alpha}{(\alpha+1)\mu^{\alpha+1}} t^{\alpha+1}\right)\\=&\exp\left(-\frac{1}{\mu}\int_0^t s dL_s - \frac{C_\alpha}{(\alpha+1)\mu^{\alpha+1}} t^{\alpha+1}\right),\end{align*}
and the right-hand side equals $\phi(t)$, which has mean $1$ by \cite[Proposition 3.2]{Conchon2018}. Also, $\E[\phi(n,t)]=1$. Then, by \cite[Lemma 4.8]{Conchon2018}, we must have that
$$\phi(n,t)\overset{d}{\to} \phi(t),$$
and  $\{\phi(n,t),n\in\N\}$ is uniformly integrable.\end{proof}

\begin{proof}[Proof of Proposition \ref{prop.laplaciantildeK}]
We want to show that for $\tilde{K}^\lambda_t$ as defined in \eqref{def.tildeK} by for $t\geq 0$, for every non-negative integrable functional $F:\D([0,t],\R^2)\to\R$, 
$$\E[F((\Tilde{K}^{\lambda}_s), 0\leq s \leq t)]=\E\left[\phi(t)F((L_s+\lambda s),0\leq s \leq t)\right],$$ 
we have that
$$\E\left[\exp\left(-\theta \tilde{K}^\lambda_t \right)\right]=\exp\left( -\theta \lambda t + \theta C_\alpha \frac{t^\alpha}{\mu^\alpha}+\int_0^t ds \int_0 ^\infty \frac{c}{\mu}x^{-(\alpha+1)}e^{-xs/\mu}dx(e^{-\theta x}-1+\theta x)\right) .$$
The proof will be an adaptation of the proof of \cite[Proposition 6.2]{Conchon2018}.
As before, let $L$ be the spectrally positive $\alpha$-stable L\'evy process having L\'evy measure $\pi(dx)=\frac{c}{\mu}x^{-(\alpha+1)}dx$ and Laplace transform
$$\Upsilon(\theta)=\frac{C_\alpha}{\mu}\theta^\alpha.$$
Let $X_t$ be the process with Laplace transform
$$\E\left[\exp(-\theta X_t)\right]=\exp\left(\int_0^t ds \int_0^\infty \pi(dx)\left(e^{-\theta x}-1+\theta x\right) e^{-\theta x s}\right),$$
and let 
$$A_t=-C_\alpha \frac{t^\alpha}{\mu^\alpha}= -\mu\Upsilon\left(\frac{t}{\mu}\right).$$
Set $$\tilde{L}^\lambda_t=X_t+A_t+\lambda t.$$
We claim that 
$$(\tilde{L}^\lambda_t,t\geq 0)\overset{d}{=}(\tilde{K}^\lambda_t,t\geq 0).$$
Like in \cite{Conchon2018}, decomposing $\pi$ gives that
$$\E[\exp(-\theta X_t)]=\exp\left(-\mu \theta \Upsilon\left(\frac{t}{\mu}\right)+\int_0^t\Upsilon\left(\theta+\frac{s}{\mu}\right)ds-\int_0^t \Upsilon\left(\frac{s}{\mu}\right)ds\right)$$
such that
\begin{align*}\E[\exp(-\theta \tilde{L}^\lambda_t)]&=\exp\left(\int_0^t\Upsilon\left(\theta+\frac{s}{\mu}\right)ds-\int_0^t \Upsilon\left(\frac{s}{\mu}\right)ds -\theta\frac{t}{\mu}\right)\\
&=\E\left[\exp\left(-\int_0^t \left(\theta+\frac{s}{\mu}\right)dL_s -\int_0^t\Upsilon\left(\frac{s}{\mu}\right)ds -\theta\frac{t}{\mu}\right)\right].\end{align*}
Let $0=t_0<t_1<\cdots<t_m=t$ and let $\theta_1,\dots,\theta_m\in \R_+$. Then, since $X$ has independent increments, by the argument as presented above,
\begin{align*}
&\E\left[\exp\left(-\sum\limits_{i=1}^m \theta_i\left(\tilde{L}^\lambda_{t_i}-\tilde{L}^\lambda_{t_{i-1}}\right)\right)\right]\\&=
\prod\limits_{i=1}^m\E\left[\exp\left(-\int_{t_{i-1}}^{t_i}\left(\theta_i+\frac{s}{\mu}\right)dL_s -\int^{t_i}_{t_{i-1}}\Upsilon\left(\frac{s}{\mu}\right)ds -\frac{\theta_i}{\mu}(t_i-t_{i-1})\right)\right]\\
&=\E\left[\exp\left(-\sum\limits_{i=1}^m \int_{t_{i-1}}^{t_i} \left(\theta_i+\frac{s}{\mu}\right)dL_s -\int^{t}_{0}\Upsilon\left(\frac{s}{\mu}\right)ds -\sum\limits_{i=1}^m \frac{\theta_i}{\mu}(t_i-t_{i-1})\right)\right] \\
&=\E\left[\exp\left(-\sum\limits_{i=1}^m \theta_i \left((L_{t_i}+\lambda t_i)-(L_{t_{i-1}}+\lambda t_{i-1})\right) -\frac{1}{\mu}\int_0^t s dL_s - \int_0^t \Upsilon\left(\frac{s}{\mu}\right)\right)\right],
\end{align*}
where we use that $L$ also has independent increments and integration by parts. This proves the statement.\end{proof}

\bibliographystyle{acm}
\bibliography{Bibliography.bib}

\end{document}